\documentclass[preprint,12pt]{elsarticle}

\usepackage{amssymb}
\usepackage{amsmath}
\usepackage{graphicx}%
\usepackage{hyperref}
\usepackage{algorithm}%
\graphicspath{{figs/}}
\usepackage{algorithmicx}%
\usepackage{algpseudocode}%
\usepackage{enumitem}
\usepackage{xcolor}%
\usepackage[section]{placeins}
\usepackage{natbib}

\journal{Elsevier}

\begin{document}

\begin{frontmatter}



\title{Quantifying Non-linearity in Topology Optimization with similarity based Visualization}


\author[1]{Ziliang Wang}
\author[1]{Jiahua Wu}
\author[1]{Jun Yang}

\author[1]{Shintaro Yamasaki\corref{cor1}}
\ead{s_yamasaki@waseda.jp}

\affiliation[1]{organization={Graduate School of Information, Production and Systems, Waseda University},
            addressline={2-7 Hibikino, Wakamatsu}, 
            city={Kitakyushu},
            postcode={808-0135}, 
            state={Fukuoka},
            country={JAPAN}}

\cortext[cor1]{Corresponding author.}

\begin{abstract}

Topology optimization (TO) can be viewed as seeking an optimal solution in the design space of a given TO problem. For weakly non-linear TO problems, e.g., compliance minimization, sensitivity-based methods typically converge well, whereas for strongly non-linear problems, e.g., maximum stress minimization, stabilization strategies such as stabilization terms and projection functions are often required to enhance convergence. Especially in scenarios with massive design variables, it is difficult to intuitively demonstrate the non-linear complexity of different TO problems and to elucidate the mechanisms by which stabilization strategies affect convergence. To address this challenge, we propose a visualization framework and a quantitative non-linearity index for objectives with varying complexity. We employ a multi-start fixed-gradient sampling tailored to similarity-based dimensionality reduction while keeping the computational cost under control. The samples are then parameterized via cosine similarity to obtain a low-dimensional visualization surface of the objective function. Based on this visualization, we construct a dimensionless complexity index with a clear geometric interpretation by measuring the gap between the visualization surface and the discrete approximation of its convex envelope, which enables quantitative comparisons of non-linearity across TO tasks, parameter choices, and stabilization strategies. Extensive comparative experiments show that the proposed approach is both adaptable and discriminative on a variety of representative TO problems, and it provides intuitive and measurable guidance for parameter selection.

\end{abstract}

\begin{graphicalabstract}

\centering
\includegraphics[width=0.9\textwidth, bb=426 150 3500 3882]{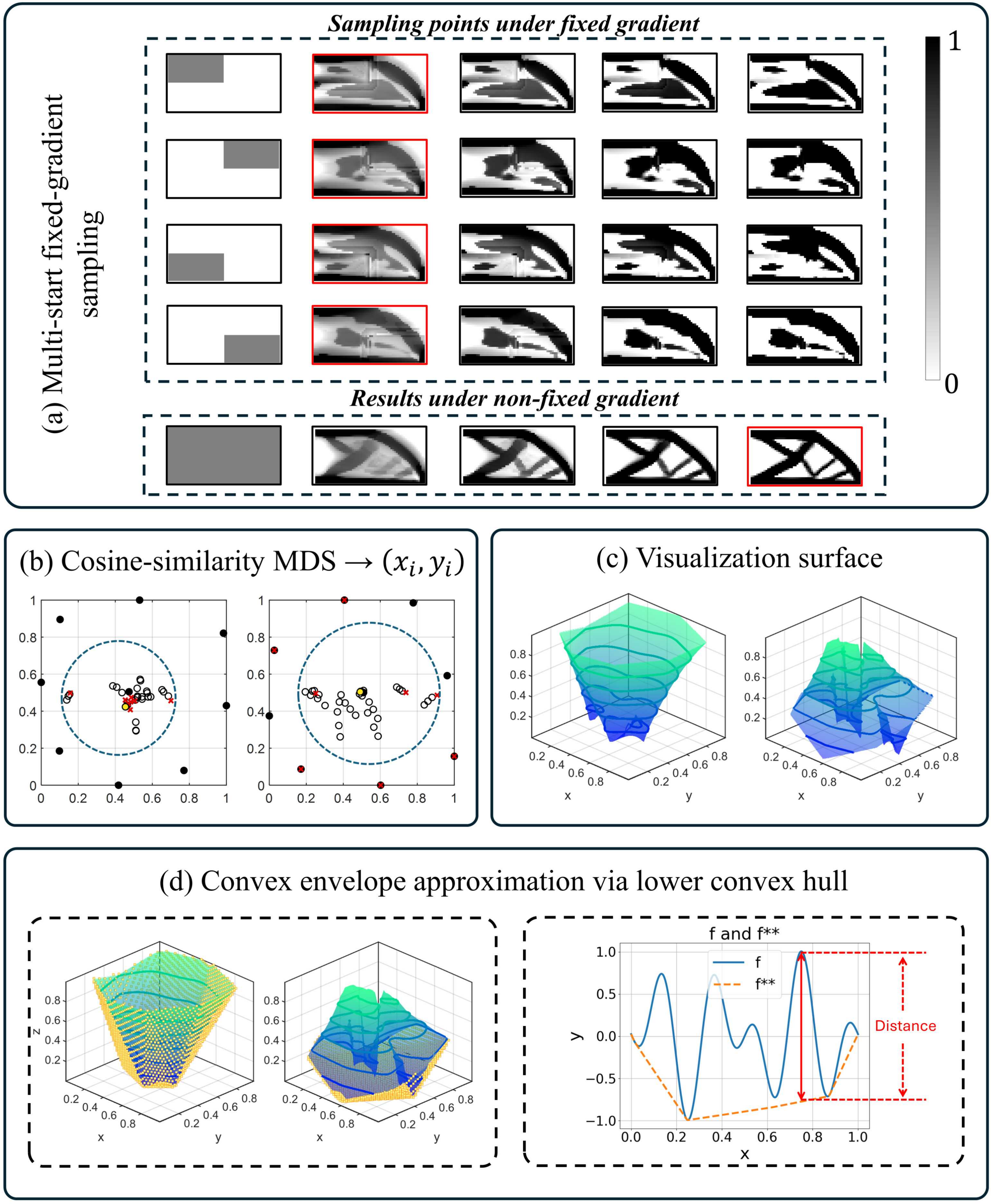}

\end{graphicalabstract}

\begin{highlights}
\item Proposed similarity-based visualization for topology optimization problems.

\item Developed a dimensionless non-linearity index via distance to lower convex hull.

\item Analyzed the effect of constraints on topology optimization non-linear complexity

\item Explained stabilization and guided parameter selection via the index.

\item Validated method effectiveness through extensive multiphysics TO benchmarks.

\end{highlights}


\end{frontmatter}



\begin{figure*}[!htp]
\centering
\includegraphics[width=0.8\textwidth, keepaspectratio, bb=0 0 3869 1298]{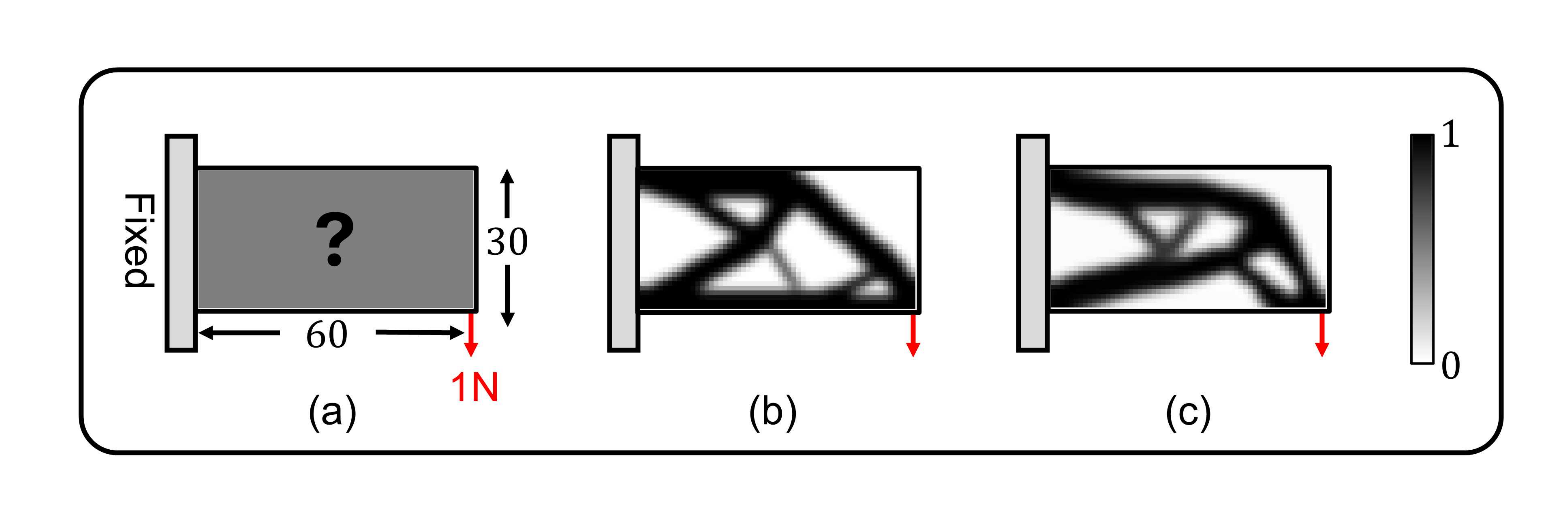}

\caption{Problem setting and optimal material distribution. (a) Design domain (dark gray) discretized by $30\times60$ elements; the left boundary is fixed and a downward distributed load of $1\,\mathrm{N}$ is applied at the lower-right corner. (b) Optimal material distribution for compliance minimization. (c) Optimal material distribution for the $p$-norm approximation of the von Mises stress objective with $p=10$. Optimization settings: optimizer, Method of Moving Asymptotes (MMA); initial element density, $0.5$; maximum allowable volume fraction, $0.5$; maximum iterations, $100$; density-filter radius, $2.5$ (element size). Grayscale shows density (white $=0$, black $=1$). 
}
\label{setup}
\end{figure*}

\section{Introduction}
\label{sec1}

Topology optimization (TO) seeks an optimal material distribution within a given design domain to minimize objectives subject to constraints. Unlike shape or size optimization, where the design variables are boundary coordinates or size variables (e.g., cross-sectional areas or shell thicknesses), TO formulations permit topological changes, allowing holes to appear or disappear during the optimization process, thereby exploring high-performance structures that are difficult to obtain through conventional empirical design \cite{eariest_TO,redbook,ole_review}. Among various TO methods, density-based approaches \cite{densitystart} are widely used because they are simple to implement and numerically reliable, and the number of design variables typically does not change during the optimization. In discretized formulations, one density variable is assigned to each finite mesh element, so the design dimension equals the number of mesh cells and grows with mesh resolution, which yields a high-dimensional design space. To enable sensitivity analysis, the design variables are relaxed from binary to continuous values from 0 to 1. Although intermediate densities may be meaningful in particular formulations, the final layout is expected to be nearly binary. As noted by Olesen et~al.\,\cite{ole_filter}, keeping gray elements in the final design is undesirable.
In common practice, standard penalization encourages binarization, for example the Solid Isotropic Material with Penalization (SIMP) for constitutive property interpolation in solids \cite{ole_simp,namesimp} and Brinkman penalties that increase the inverse permeability in fluid formulations \cite{dasistart}. As an illustrative example, we consider the 2D cantilever configuration in Fig.~\ref{setup}(a); Fig.~\ref{setup}(b) and (c) show material distribution optimized for two objectives, compliance and maximum stress.

\begin{figure*}[!htp]
\centering
\includegraphics[width=0.8\textwidth,keepaspectratio, bb=281 210 3944 2356]{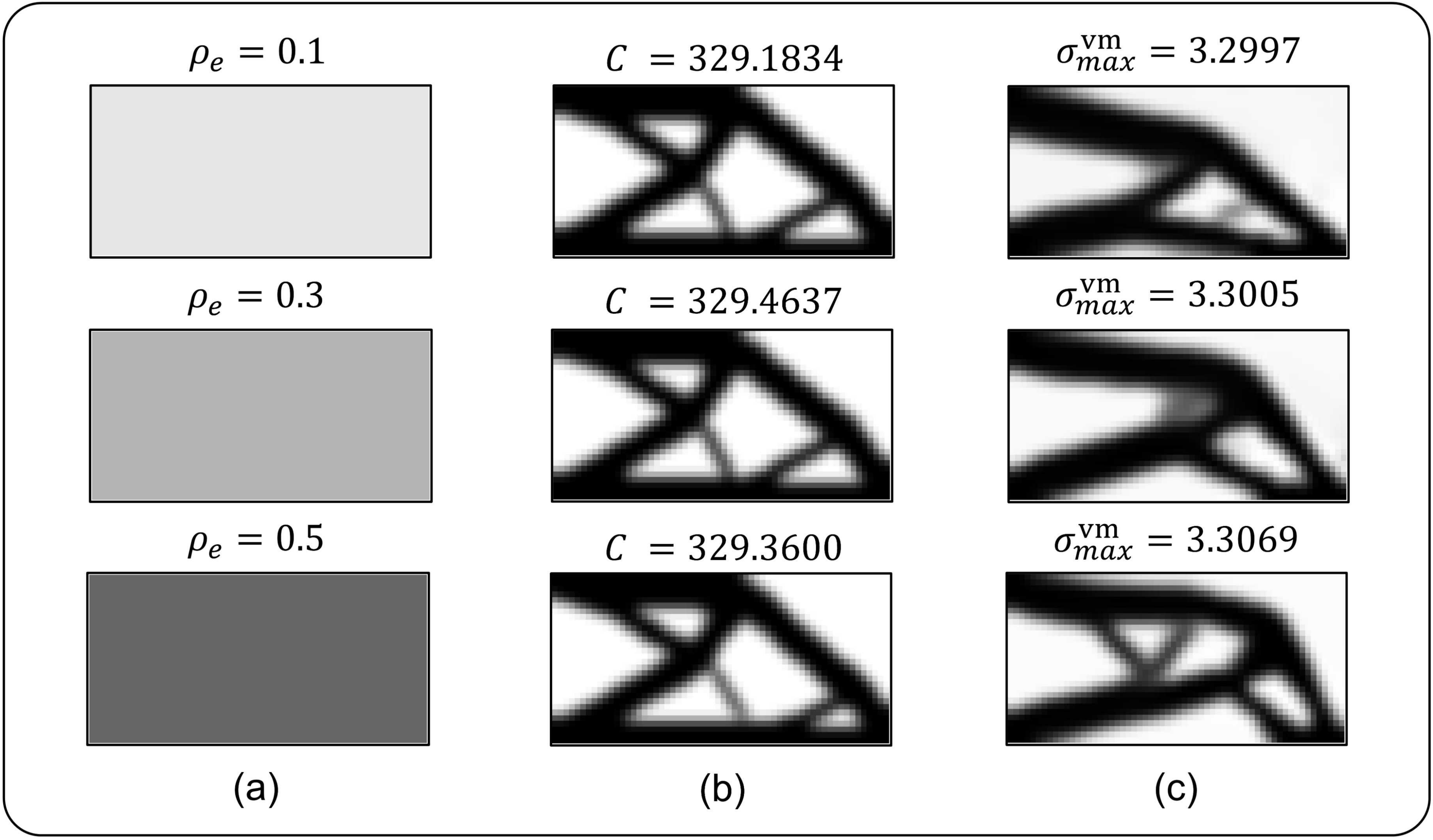}
\caption{Weakly vs.\ strongly non-linear objectives. The problem setting follows Fig.~\ref{setup}(a) and is discretized by $30\times60$ elements. Optimizer: MMA; maximum iterations: 100; density-filter radius: 2.5 (elements); $p$-norm stress with $p=10$. (a) Three uniform initial designs with element density $\rho_e \in \{0.1,0.3,0.5\}$. (b) Compliance minimization, started from (a), converges to similar local optima. (c) Maximum von Mises stress minimization, due to higher non-linearity, converges to distinct local optima.}

\label{distinct_local_optima}
\end{figure*}

For weakly non-linear TO problems, such as compliance minimization, convergence is usually robust, i.e., different uniform initial densities in Fig.~\ref{distinct_local_optima}(a) lead to similar local optima in Fig.~\ref{distinct_local_optima}(b). However, for highly non-linear TO problems, e.g., maximum stress minimization, the optimization process often exhibits slow convergence and numerical instability \cite{stressconvergeslow}. Even under identical constraints and design domain, different initializations can derive obviously different local optima, as shown in Fig.~\ref{distinct_local_optima}(c). To further promote convergence and mitigate such instabilities, prior studies have proposed strategies such as introducing stabilization terms \cite{3field2024}, using projection functions \cite{3field2022}, and heuristic parameter tuning \cite{3field2025}, but the mechanisms behind these strategies still lack intuitive explanations and quantitative evaluation frameworks.

Meanwhile, density-based methods involve massive parameters. Taking maximum stress minimization as an example, typical parameters include the $p$ in the $p$-norm aggregation, the filter radius, the Heaviside projection parameters, and the penalty and relaxation coefficients. For such problems are highly non-linear and often ill-conditioned, the optimized designs are also strongly dependent on mesh resolution \cite{3field2024}. As the mesh resolution increases, different discretizations may converge to different local optima, while the number of design variables and the simulation cost both grow substantially \cite{olecostwithmesh}. Consequently, trial-and-error parameter selection becomes computationally expensive. Also, designs obtained on coarse meshes cannot be scaled directly to finer meshes \cite{olegooddesignscale}. To alleviate these issues, the Three-field floating projection topology optimization (FPTO) has been proposed \cite{3field2024}. By heuristically adjusting the Heaviside projection parameters and introducing the average solid stress, it reduces sensitivity to the choice of $p$ and to mesh resolution, thereby improving convergence speed and robustness. Nevertheless, selecting a suitable combination of parameters still relies heavily on experience and lacks intuitive, quantitative selection criteria.

Visualization offers a promising avenue to address the above challenge. If the landscape of the objective function can be visualized effectively, the practical impact of different stabilization strategies and parameter combinations on problem complexity can be analyzed, thereby reducing evaluation costs \cite{qualitymetricsvisual,design_space_dimensionality_reduction}. However, in density-based TO, the optimization is performed in a high-dimensional design space, whose dimension equals the number of design variables. Direct analysis and comparison in this space are both computationally expensive and difficult to interpret. Consequently, constructing low dimensional representations of the design space has become an active research direction. For instance, incorporating principal component analysis (PCA) \cite{pca} into deep learning based TO can alleviate the inherent input-dimensionality bottleneck of deep generative models \cite{PCADDTD}. Another research construct mappings from the high dimensional finite element space to a low dimensional manifold, compute convex hulls in the low dimensional space, perform resampling, and then inverse map to the high dimensional space to generate new candidates and accelerate optimization \cite{PCAacceleratingphase-field}. Moreover, combining manifold learning with deep neural networks allows training predictors on the low dimensional manifold for rapid reconstruction of high dimensional finite element solutions, showing efficiency in inverse problems \cite{reduceDiminverse}. Nevertheless, merely constructing a low dimensional embedding is insufficient for effective visualization and analysis of the objective function, i.e., if the embedding does not preserve relative positional relationships from the original space and only displays scalar function values, it becomes difficult to reliably quantify and compare complexity across problems, stabilization strategies, and parameter settings.

In deep learning, visualizing the loss landscape of neural networks has proved valuable. Studies report that, in latent spaces, there are low-dimensional linear subspaces along which the loss appears nearly convex \cite{goodfellowlandscape}. Other work compares how network architectures affect loss landscape complexity \cite{landscapevisualizing}. These studies demonstrate the potential of visualization for exploring loss function complexity. However, sampling along straight lines between two parameter points rarely follows the trajectories taken by an optimizer, so its ability to reflect optimization difficulty is limited. In addition, these visualizations typically ignore constraints, since standard network training is formulated as an unconstrained problem, and they do not provide a dimensionless non-linearity measure that enables systematic comparisons across tasks and parameter settings.
In contrast, quantitative complexity analysis for TO problems must account for both the objective and the constraint functions, such as volume. Therefore, effective visualization thus relies on representative samples that satisfy the relevant constraints, capture the geometric structure of high-dimensional objective functions with manageable computational cost, and preserve key relative positional information under aggressive dimensionality reduction. 

From optimization theory, non-linearity is closely linked to problem complexity. 
For convex optimization problems, any local optimum is also global, and many standard algorithms enjoy well established convergence guarantees. 
In contrast, highly non-linear TO objectives typically admit multiple local optima, so the gap between the objective and its tight convex envelope provides a natural measure of non-linearity. 
The Legendre-Fenchel transform (LFT) is a classical tool, and applying it twice yields the tightest convex envelope of a function \cite{LFT1,LFT2}. 
However, for high-dimensional TO objectives, explicitly computing the LFT biconjugate via two successive transforms is prohibitively expensive.

\begin{figure*}[!htp]
\centering
\includegraphics[width=0.8\textwidth, keepaspectratio, bb=215 181 3857 1575]{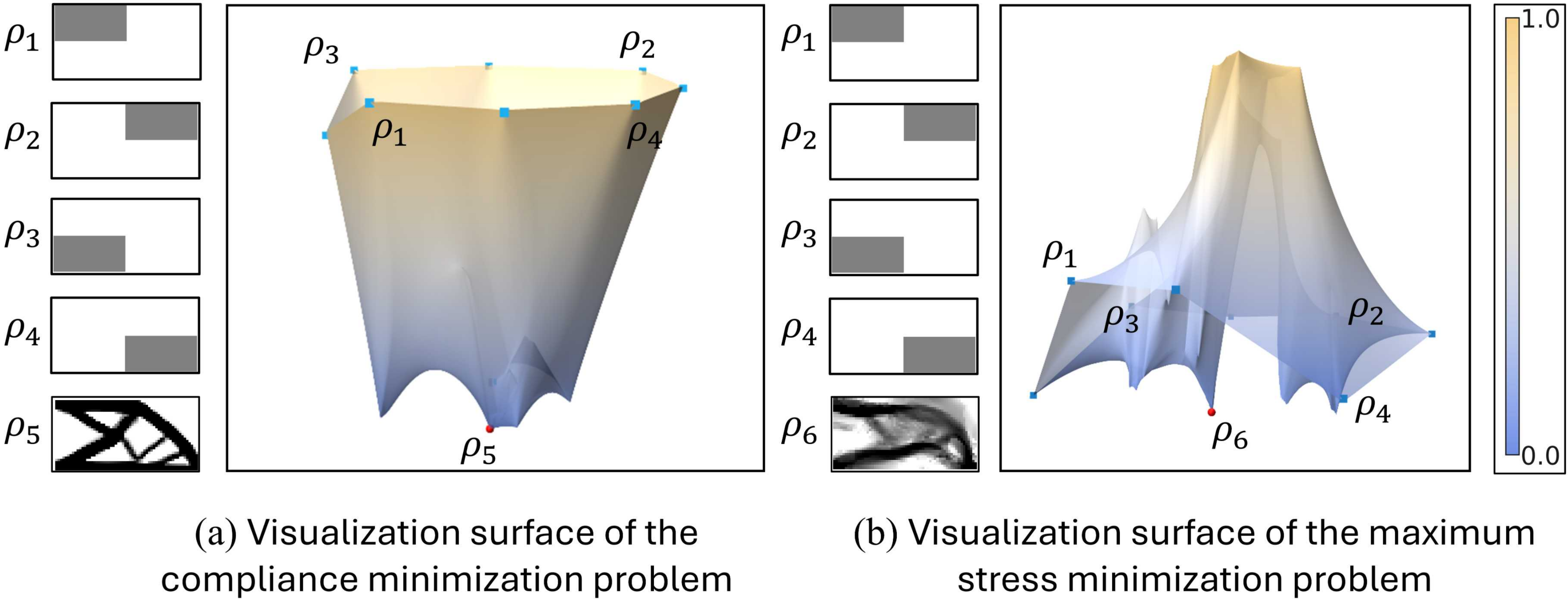}
\caption{Comparison of visualization surfaces. The problem setting follows Fig.~\ref{setup}(a), except that the move limit control parameter $\eta_{\max}$ is reduced to $0.01$ after the fifth iteration, which leads to minor differences in the optimal material distribution and this modification is discussed in Section~\ref{sec:ablationaaa}. To preserve the cosine distances between initial points, we use the unfiltered design variables as sampling points. Left, compliance as the objective yields a relatively simple surface with few peaks and valleys. Right, the $p$-norm approximation of the maximum von Mises stress objective ($p=10$) is visibly more complex. Blue squares indicate the starting points for each sampling group, with the red dot indicating the best performing design.}

\label{visualization}
\end{figure*}

On the basis of the above discussions, we develop a visualization and quantitative evaluation framework for constrained TO objectives that addresses three practical challenges. (see Fig.~\ref{visualization} for examples of the visualization results). To obtain representative samples that are compatible with the constraints and reflect TO problem complexity, we employ multi-start fixed-gradient sampling to generate task-aware samples under a controlled computational budget. To make visualization and comparison reliable in spite of aggressive dimensionality reduction, we construct a multidimensional scaling (MDS) \cite{mds} embedding using cosine similarity distances, which preserves relative positional relationships among designs from the original high dimensional space. To quantify non-linearity at controlled computational cost, we define a dimensionless index based on the gap between the visualization surface and its convex envelope, where the convex envelope is approximated by the discrete lower convex hull, avoiding the prohibitive cost of two successive Legendre-Fenchel transforms in high dimensions.

The main contributions are as follows.
First, we propose a visualization framework for constrained TO objectives that combines multi-start fixed-gradient sampling with an MDS embedding based on cosine similarity, yielding high fidelity visualization surfaces at a controlled computational cost.
Second, we construct a dimensionless non-linearity index from the gap between the visualization surface and its convex envelope, which is approximated by the discrete lower convex hull.
Third, systematic studies on structural, heat transfer, and fluid problems demonstrate that the proposed framework and index distinguish differences in non-linearity across objectives and stabilization strategies, are robust over a range of parameter settings, and support intuitive and quantitative interpretation of stabilization effects and parameter selection.

The remainder of this paper is organized as follows. Sect. \ref{sec:methods} details the proposed method, Sect. \ref{section3} presents numerical experiments, Sect. \ref{sec:restrictions_outlook} discusses conclusion and limitations.

\section{Methods}
\label{sec:methods}

We consider a density-based TO problem on a design domain discretized into \(n\) finite elements. The design vector is \(\boldsymbol{\rho}\in\mathbb{R}^n\) and it satisfies the constraints \(\rho_{\min}\le \rho_i\le 1\) for all \(i\), where \(\rho_{\min}=\varepsilon\) is a small positive regularization parameter (e.g., \(\varepsilon\approx 10^{-6}\)). Let \(\mathbf{s}\) denote the state variables of the governing physics, e.g., displacement \(\mathbf{u}\) in structures, temperature \(T\) in heat conduction, and velocity \(\mathbf{v}\) and pressure \(p\) in fluid problems. The general form of the TO problem is formulated by

\begin{equation}
\label{eq:TO}
\begin{aligned}
\min_{\boldsymbol{\rho}\in[\rho_{\min},1]^n} \quad 
& J\big(\mathbf{s}(\boldsymbol{\rho}),\,\boldsymbol{\rho}\big) \\[2pt]
\text{s.\,t.}\quad 
& \mathbf{R}\big(\mathbf{s}(\boldsymbol{\rho}),\,\boldsymbol{\rho}\big)=\mathbf{0}, \\[2pt]
& g_k(\boldsymbol{\rho}) \le \bar g_k, \qquad k=1,\dots,m .
\end{aligned}
\end{equation}

Here \(\mathbf{R}(\cdot,\cdot)=\mathbf{0}\) denotes the discrete residual of the governing equations, and \(J\) is the objective function. As an illustration in linear elasticity, $\mathbf{s}=\mathbf{u}$ and $\mathbf{K}(\boldsymbol{\rho})\,\mathbf{u}=\mathbf{f}$, while the compliance objective is $J(\mathbf{u},\boldsymbol{\rho})=\mathbf{f}^{\top}\mathbf{u}$, 
where $\mathbf{K}(\boldsymbol{\rho})\in\mathbb{R}^{n_{\mathrm{dof}}\times n_{\mathrm{dof}}}$ is the global stiffness matrix and $\mathbf{f}\in\mathbb{R}^{n_{\mathrm{dof}}}$ is the global load vector, with $n_{\mathrm{dof}}$ denoting the number of degrees of freedom.

\begin{figure*}[!htp]
\centering
\includegraphics[width=0.8\textwidth, keepaspectratio, bb=426 150 3500 3882]{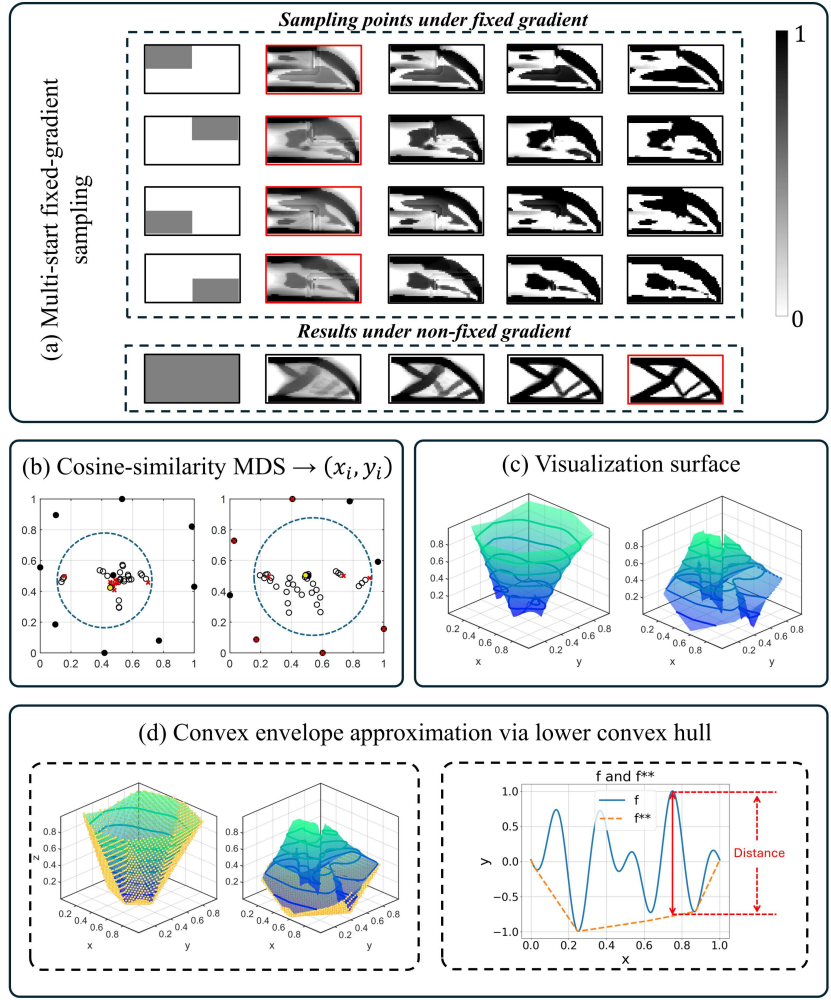}
\caption{Schematic of the proposed framework. (a) Multi-start fixed-gradient sampling under a fixed gradient from different initial points; within each group, the best result is highlighted with a red box. In (b), (c), and the left of (d), the left column shows the compliance minimization problem, and the right column shows the minimization of the maximum von Mises stress. (b) Cosine-based MDS maps the high-dimensional designs to 2D coordinates. (c) Adding normalized objective values on the embedding yields the visualization surface. (d) Left: visualization surface with sample points on its convex envelope. Right: a 2D toy example illustrating the distances between the visualization surface and its convex envelope.
}
\label{pipeline}
\end{figure*}

For density-based TO problems, our goal is to visualize the landscape of the objective function subject to constraints and use it to define a quantitative, dimensionless index of non-linearity. Fig.~\ref{pipeline} outlines the workflow: (a) Multi-start fixed-gradient sampling collects representative designs at controlled computational cost, yielding samples that reveal and differentiate problem complexity;
 (b) the designs are embedded into a low dimensional space using classical MDS based on cosine similarity; (c) We assign normalized objective function values to the embedding coordinates to generate a visualization surface; and (d) shows a discrete approximation of the convex envelope of the surface, obtained from the lower convex hull, where vertical distance between the surface and this envelope define the non-linearity index.
 The remainder of this section details each stage in turn.

\subsection{Multi-start fixed-gradient sampling}
\label{subsec:fixed_gradient}

To quantify problem complexity more accurately, we require a sampling scheme that \emph{reveals non-linearity}. In high-dimensional design spaces, purely random sampling offers little control over coverage, and uniform (grid) sampling grows rapidly with the dimension, so achieving accuracy typically requires a large number of samples and substantially increases the cost of subsequent computations.

\begin{figure}[!htp]
\centering
\includegraphics[width=0.8\textwidth, keepaspectratio, bb= 254 228 3100 3951]{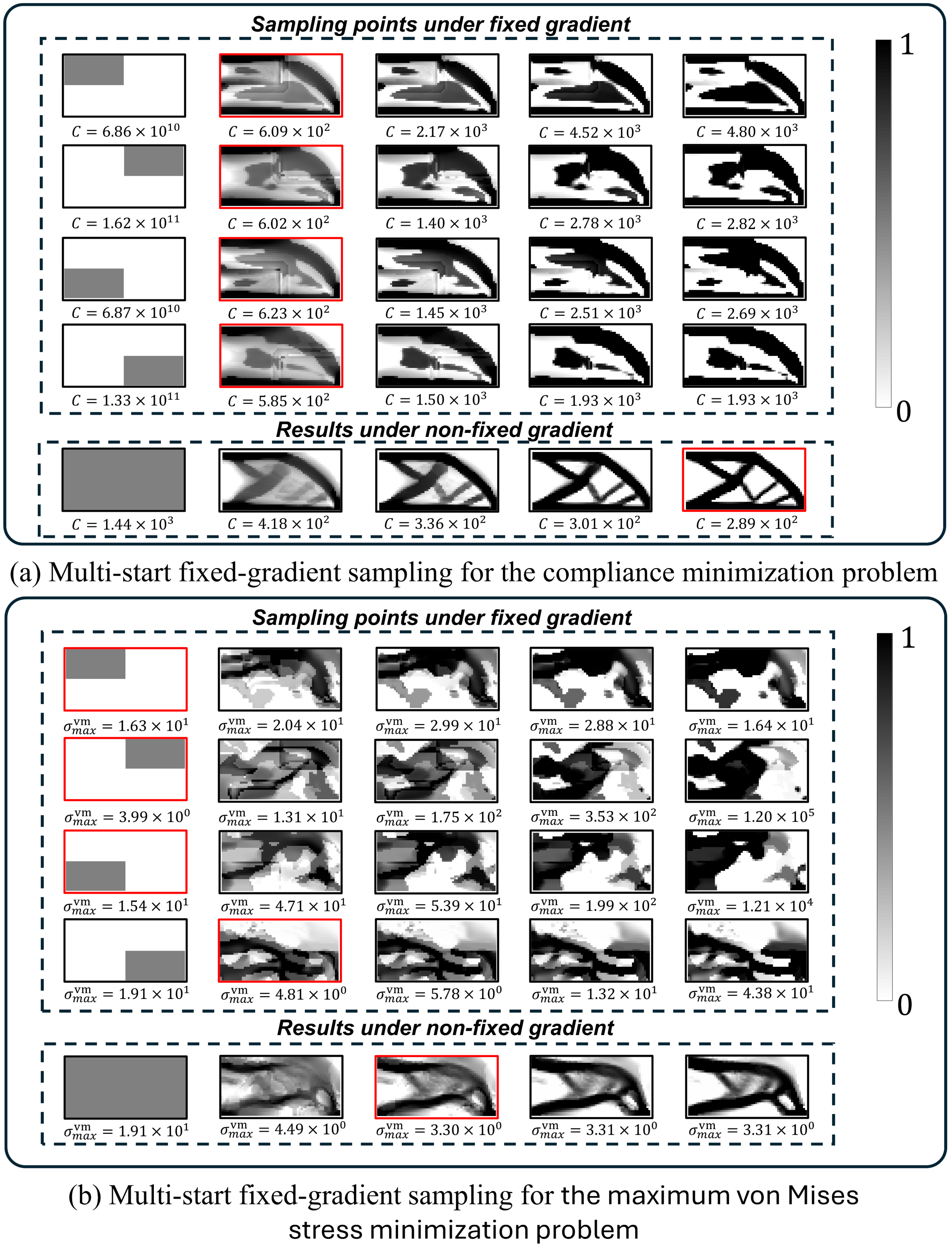}
 
\caption{Multi-start fixed-gradient sampling. Starting from different initial points, we iterate with a fixed gradient for 100 steps and save the iterate every 20 steps. In both the top and bottom panels, each row corresponds to one group, with the red box marking the best structure in that group. We observe that, for (a), the per-group best structures are more similar to the final designs obtained by standard optimization without the fixed-gradient constraint (bottom strip), whereas for (b), the best structures differ markedly from those final designs in shape. These differences also motivate our use of a similarity-based dimensionality reduction method in the subsequent stage. The problem setting is identical to Fig.~\ref{visualization}.}

\label{sample}
\end{figure}

We consider $M$ starting designs 
$\{\boldsymbol{\rho}^{(j)}_0\}_{j=1}^M \subset [\rho_{\min},1]^n$.
For each start, we first run a standard topology optimization (TO) optimizer
for $t_{\mathrm{fix}} \ge 0$ iterations to mitigate gradient instability in the early optimization stages.
The objective gradient is then evaluated at
$\boldsymbol{\rho}^{(j)}_{t_{\mathrm{fix}}}$ and kept fixed,
\[
  \mathbf{J}^{(j)}_{\mathrm{fix}}
  =
  \nabla_{\boldsymbol{\rho}} J\!\big(
  \mathbf{s}(\boldsymbol{\rho}^{(j)}_{t_{\mathrm{fix}}}),
  \boldsymbol{\rho}^{(j)}_{t_{\mathrm{fix}}}
  \big).
\]
In all numerical experiments we set $t_{\mathrm{fix}} = 5$.

We define the projection $\Pi:\mathbb{R}^n \to [\rho_{\min},1]^n$
onto the box constrained design space in a component-wise form,
\[
  \bigl(\Pi(x)\bigr)_i
  =
  \min\!\bigl(\max(x_i,\rho_{\min}),\,1\bigr),
  \qquad i=1,\ldots,n.
\]

Let $\mathcal{A}_{\mathrm{fg}}(\cdot;\mathbf{J}^{(j)}_{\mathrm{fix}},\eta_{\max})$
denote one update step of a gradient-based constrained optimizer in which
the \emph{objective gradient} is fixed at $\mathbf{J}^{(j)}_{\mathrm{fix}}$,
while the objective value, constraint values, and constraint gradients
are still evaluated at the current iterate.
The parameter $\eta_{\max}$ corresponds to the \emph{move limit} used in the optimizer,
for example the \texttt{move} parameter in the Method of Moving Asymptotes (MMA),
which controls the maximum allowable change of design variables per iteration; see \cite{svanberg2007mma}.

Given a total number of iterations $T \in \mathbb{N}$,
the fixed-gradient sampling path is defined as
\[
  \boldsymbol{\rho}^{(j)}_{t}
  =
  \Pi\!\big(
  \mathcal{A}_{\mathrm{fg}}(
  \boldsymbol{\rho}^{(j)}_{t-1};
  \mathbf{J}^{(j)}_{\mathrm{fix}},
  \eta_{\max})
  \big),
  \qquad t = t_{\mathrm{fix}}+1,\ldots,T.
\]

Let $\Delta \in \mathbb{N}$ denote the save interval.
Along each sampling path we store the designs
\[
  \big\{
  \boldsymbol{\rho}^{(j)}_{q\Delta}
  \;\big|\;
  q = 0,1,\ldots,\big\lfloor T/\Delta \big\rfloor
  \big\},
\]
and evaluate the corresponding objective values
\[
  J^{(j)}_{q\Delta}
  =
  J\!\big(
  \mathbf{s}(\boldsymbol{\rho}^{(j)}_{q\Delta}),
  \boldsymbol{\rho}^{(j)}_{q\Delta}
  \big),
  \qquad
  q = 0,1,\ldots,\big\lfloor T/\Delta \big\rfloor,
\]
where $\lfloor\cdot\rfloor$ denotes the floor operator. Fig.~\ref{sample}(a) and (b) present the sampling points of each group for the compliance minimization problem and the maximum von Mises stress minimization problem, respectively.

Having introduced the sampling method above, two design choices warrant clarification:
(i) why sampling is performed via a constrained optimizer, and
(ii) why the objective gradient is fixed after $t_{\mathrm{fix}}$ iterations.

\paragraph{(i) Why a constrained optimizer}
TO problems differ from the loss landscape visualization commonly discussed in deep learning because explicit constraints (for example, a volume constraint such as keeping the material volume fraction below $50\%$) must be respected for sampled material distributions to remain task-aware.
Consequently, straightforward linear interpolation between two parameter matrices, which in TO problems amounts to interpolating between two material distributions, can generate samples that violate the imposed constraints, therefore not comparable to samples produced under the TO formulation.
We therefore perform sampling via a constrained optimizer.
Here, the optimizer is not used to seek an optimal solution.
Instead, it offers a simple and practical way to take constraints into account during sampling, so that the resulting samples remain aligned with the underlying TO task.
This is essential because our goal is to reveal how constraints influence the non-linearity of TO problems, rather than to generate interpolated samples that disregard the constraints.

\paragraph{(ii) Why the gradient is fixed}
Fixing the objective gradient provides two practical benefits.
First, it leads to more dispersed sampled designs and reduces sample overlap, which is particularly beneficial for TO problems with low non-linearity.
In such cases, sampling without fixing the objective gradient tends to produce highly concentrated samples, potentially obscuring differences in problem complexity.
Second, fixing the gradient reduces the sampling cost by avoiding repeated evaluations of the objective gradient.
The impact of fixing the objective gradient is further investigated and quantitatively demonstrated in the ablation study presented in Section~\ref{sec:ablationaaa}.

\subsection{MDS and cosine similarity}
\label{subsec:mds_cosine}

As described in Section~\ref{subsec:fixed_gradient}, multi-start fixed-gradient sampling yields
a collection of saved designs
\(
\{\boldsymbol{\rho}^{(j)}_{q\Delta}\}
\),
where $j=1,\dots,M$ indexes the starting designs and
$q = 0,1,\ldots,\lfloor T/\Delta \rfloor$
indexes the saved iterations along each sampling path.
For notational simplicity, we relabel all saved designs with a single index
$i = 1,\dots,N$,
where
\(
N = M \bigl(\lfloor T/\Delta \rfloor + 1\bigr)
\)
is the total number of samples.
Each sampled design $\boldsymbol{\rho}_i$ is associated with a corresponding objective value $J_i$.

We then embed the designs into a low-dimensional space for visualization and geometric analysis. The high-dimensional dataset is
\begin{equation}
\label{eq:sample_set}
\mathcal{X}
=
\Big\{\,(\boldsymbol{\rho}_i,\, J_i)\ \Big|\ \boldsymbol{\rho}_i\in[\rho_{\min},1]^n,\ J_i\in\mathbb{R},\ i=1,\dots,N\,\Big\}.
\end{equation}
Only the design vectors \(\{\boldsymbol{\rho}_i\}_{i=1}^N\) are used to build the distance matrix for dimensionality reduction; the embedding uses only the designs, while the objective values \(\{J_i\}_{i=1}^N\) are not used at this stage. They are incorporated in the next step to construct the visualization surface.

\paragraph{Cosine distance}
To focus on geometric and topological differences rather than volume scaling, we use the cosine distance. For two designs \(\boldsymbol{\rho}_i\) and \(\boldsymbol{\rho}_j\), define the cosine similarity
\begin{equation}
\label{eq:cosine_sim_redef}
s_{ij}
=
\frac{\boldsymbol{\rho}_i^{\top}\boldsymbol{\rho}_j}{\|\boldsymbol{\rho}_i\|_2\,\|\boldsymbol{\rho}_j\|_2}
\ \in [0,1],
\end{equation}
the corresponding cosine distance is
\begin{equation}
\label{eq:cosine_dist}
d_{ij} \;=\; 1 - s_{ij} \;\in\; [0,1].
\end{equation}

We form the squared-distance matrix
\begin{equation}
\label{eq:D_from_cos}
\mathbf{D} \;=\; \big(d_{ij}^{2}\big)_{i,j=1}^{N}.
\end{equation}
\paragraph{Classical MDS}
Classical MDS constructs an embedding \(\mathbf{Y}=[\,\boldsymbol{y}_1,\dots,\boldsymbol{y}_N\,]\in\mathbb{R}^{p\times N}\) such that the Gram matrix \(\mathbf{Y}^\top\mathbf{Y}\) approximates the doubly centered Gram matrix \(\mathbf{G}\) built from squared distances. Let
\begin{equation}
\label{eq:dist_center}
\mathbf{H}=\mathbf{I}-\tfrac{1}{N}\,\mathbf{1}\mathbf{1}^{\top},
\qquad
\mathbf{G} \;=\; -\,\tfrac{1}{2}\,\mathbf{H}\,\mathbf{D}\,\mathbf{H}.
\end{equation}
The classical MDS problem can be formulated as
\begin{equation}
\label{eq:mds_obj}
\mathbf{Y}^{\ast}
=
\arg\min_{\mathbf{Y}\in\mathbb{R}^{p\times N}}
\ \big\lVert \mathbf{Y}^{\top}\mathbf{Y}-\mathbf{G}\big\rVert_{F}^{2}.
\end{equation}
Here \(p\) is the embedding dimension (typically \(p=2\) or \(3\)), \(\mathbf{1}\in\mathbb{R}^N\) denotes the all-ones column vector, \(\mathbf{I}\) is the \(N\times N\) identity matrix, and \(\|\cdot\|_{F}\) denotes the Frobenius norm.

\begin{figure}[!htp]
\centering
\includegraphics[width=0.8\textwidth, keepaspectratio, bb=695 250 3850 2738]{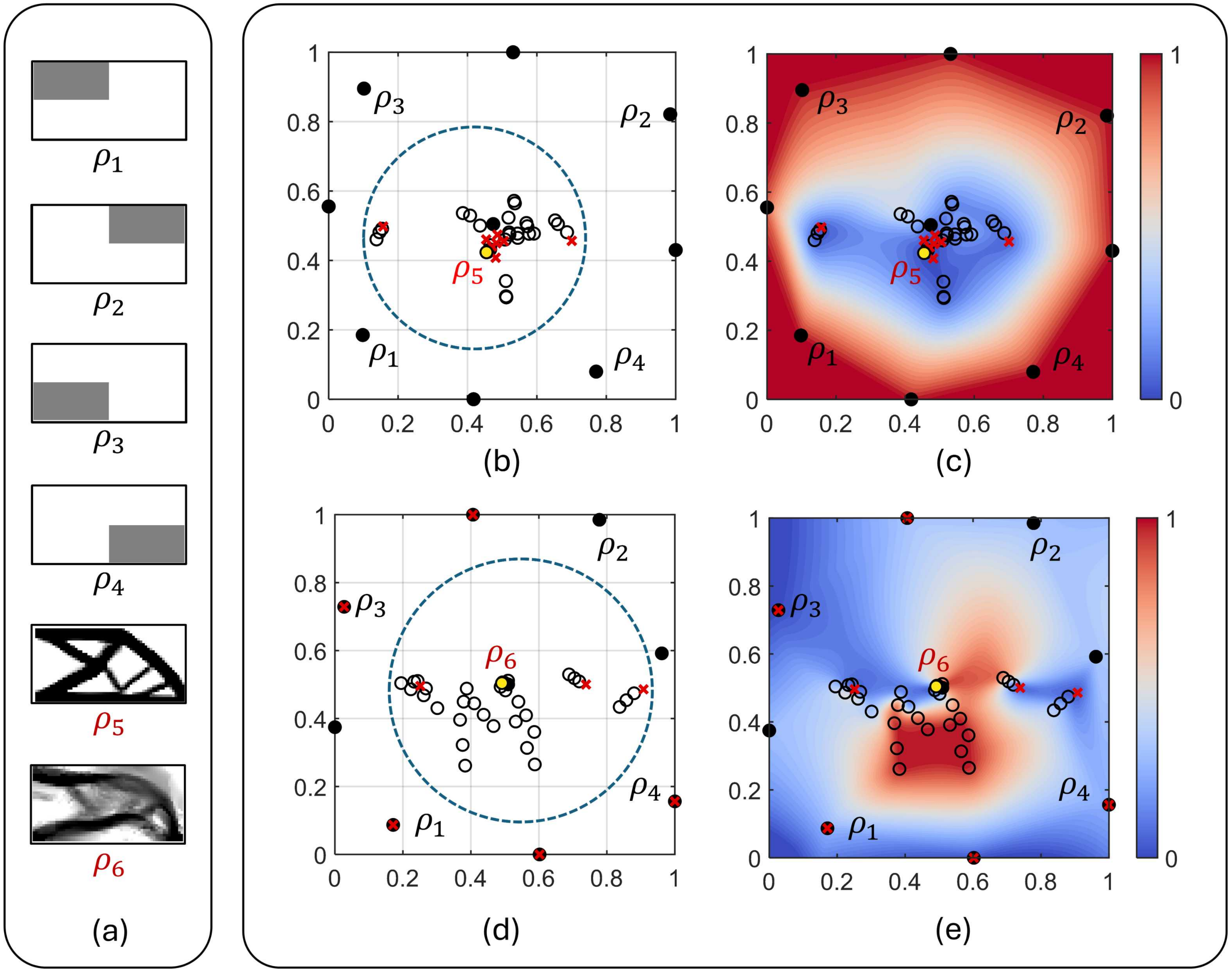}
\caption{cosine based MDS embedding of sampled designs.
(a) Selected material layouts: $\rho_{1}$ to $\rho_{4}$ are starting designs, and $\rho_{5}$, $\rho_{6}$ are reference solutions.
(b, d) Samples generated from nine different starts by multi-start fixed-gradient sampling, embedded to 2D with classical MDS using cosine distance; (b) shows the compliance minimization problem, and (d) the maximum stress minimization problem.
(c, e) The same embeddings with normalized objective values lifted onto the surface; color encodes the normalized objective value (warmer colors indicate larger values).
Each dot is a sampled design; red $\times$ mark the best sample in each group.
}
\label{cosine-distance-embedding}
\end{figure}

\paragraph{Example}
We now illustrate the above on a representative example in Fig.~\ref{cosine-distance-embedding}.
Fig.~\ref{cosine-distance-embedding}(a) lists six material layouts: from $\rho_{1}$ to $\rho_{4}$ are different starts and $\rho_{5}$ and $\rho_{6}$ are reference solutions obtained from standard topology optimization runs without fixed-gradient sampling. We include these reference solutions in the experiments for exposition, and refer to them as “reference solution” throughout. Their influence on the resulting visualization surface and on the estimated complexity is minor, as we will verify in the Section~\ref{section3}.

To keep the starts well separated under the cosine metric, we select eight pairwise orthogonal starts, i.e., every pair has cosine distance 1. To avoid distorting the constructed cosine distance relationship and to maintain consistency across all sampling points, we use the unfiltered design variables as sampling points to construct the 2D embedding.
Fig.~\ref{cosine-distance-embedding}(b) and (d) show the two-dimensional embeddings produced by classical MDS from the cosine distance on samples generated by multi-start fixed-gradient sampling; each dot denotes a sampled design, and a red $\times$ marks the best sample within each start group. Fig.~\ref{cosine-distance-embedding}(c) and (e) assign the normalized objective values on the embeddings to form an visualization surface, where color encodes the normalized objective value.
The problem setting is identical to Fig.~\ref{visualization}. Compliance minimization case in Fig.~\ref{cosine-distance-embedding}(b) and (c) yields a compact cluster, indicating a relatively low-complexity landscape. In Fig.~\ref{cosine-distance-embedding}(d) and (e), the samples are widely dispersed except for the starting points, and the reference solution is surrounded by high ridges on the visualization surface, indicating a highly non-linear problem.


\subsection{Convex envelope and lower convex hull approximation}
\label{subsec:lft}

\paragraph{Legendre-Fenchel transform (LFT)}
Let \(\varphi:\mathbb{R}^d\to\mathbb{R}\) be a not necessarily convex function.
In our setting, \(\varphi\) is the objective function of a TO problem, and \(x\in\mathbb{R}^d\) collects all design variables (e.g., element densities).
The LFT of \(\varphi\) is defined as
\begin{equation}
  \varphi^{\ast}(p)
  = \sup_{x\in\mathbb{R}^d}\{\langle x,p\rangle - \varphi(x)\},
  \qquad p\in\mathbb{R}^d,
\end{equation}
where \(\langle\cdot,\cdot\rangle\) denotes the standard Euclidean inner product.
Here \(p\) has the same dimension as \(x\) and, in our sensitivity-based TO framework, we simply take
\begin{equation}
  p = \nabla_x \varphi(x),
\end{equation}
that is, \(p\) is the gradient (sensitivity) of the objective with respect to the design variables.
The biconjugate \(\varphi^{\ast\ast}\) is obtained by applying the LFT once more to \(\varphi^{\ast}\):
\begin{equation}
  \varphi^{\ast\ast}(x)
  = \sup_{p\in\mathbb{R}^d}\{\langle x,p\rangle - \varphi^{\ast}(p)\},
  \qquad x\in\mathbb{R}^d.
\end{equation}
It is convex and satisfies
\begin{equation}
  \varphi^{\ast\ast}(x) \le \varphi(x),
  \qquad \forall\,x\in\mathbb{R}^d.
\end{equation}

\label{subsec:lf_transform_intro}

\begin{figure}[!htp]
\centering
\includegraphics[width=0.75\textwidth, keepaspectratio, bb=322 266 2720 2516]{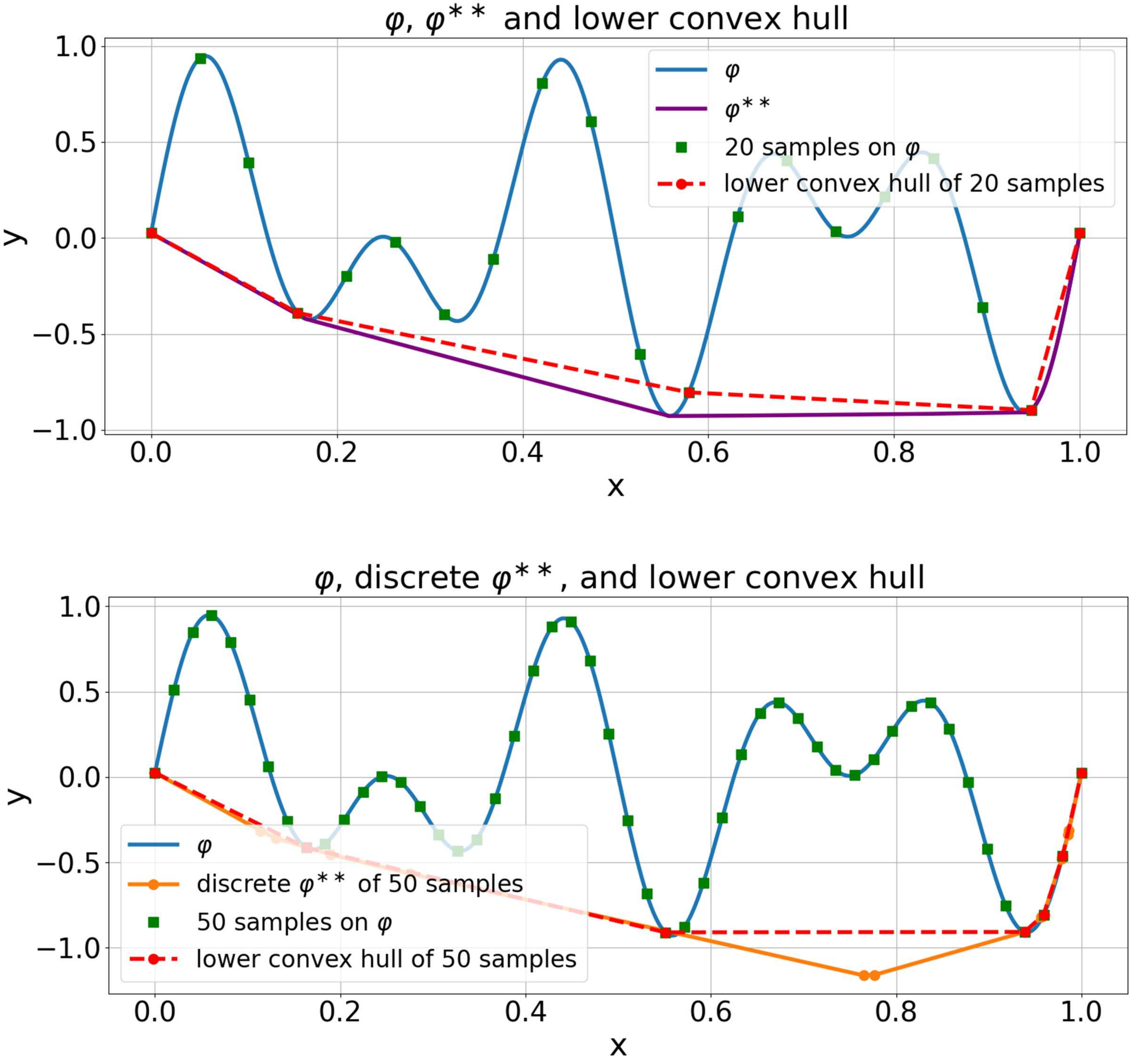}
\caption{
Illustration of a 1D non-linear function, its convex envelope, and a discrete lower convex hull. 
The blue curve shows $\varphi(x)=\sin(5\pi x)\cos(5\pi x)+\sin(3\pi x)\cos(3\pi x)+0.1(x-0.5)^2$.
The purple curve shows its convex envelope $\varphi^{\ast\ast}(x)$, computed via a high-resolution discrete double LFT.
The orange curve shows the result of applying a discrete double LFT to the original function using 50 sampled points.
The green square markers indicate uniformly sampled points on $\varphi$, and the red dashed curve with circular markers shows the lower convex hull constructed from these sampled points.}
\label{fig:lf-demo}
\end{figure}

Figure~\ref{fig:lf-demo}(a) illustrates an example in one dimension,
providing an intuitive visualization of the above concepts.
The blue curve shows a non-linear function \(\varphi\), and the purple curve shows its convex envelope \(\varphi^{\ast\ast}\), computed via a high resolution discrete double LFT.
The green markers indicate 20 samples on \(\varphi\), and the red dashed curve is the lower convex hull constructed from these points; it already closely approximates \(\varphi^{\ast\ast}\) without the two gradient evaluations per sample required by a discrete double LFT.

Figure~\ref{fig:lf-demo}(b) further compares the lower convex hull with the discrete double LFT at 50 sampled points; the lower hull tracks the convex envelope more accurately at this resolution, which motivates our use of the lower convex hull as a practical approximation of the convex envelope in our framework.

Motivated by this comparison, we approximate the convex envelope of the visualisation surface
using a lower convex hull. Visualisation surface is constructed from the sampled designs produced by
multi-start fixed-gradient sampling and cosine-based MDS (Sec.~\ref{subsec:mds_cosine}).
Each sampled design $\boldsymbol{\rho}_i$ is mapped to a two-dimensional coordinate
$\boldsymbol y_i \in \mathbb{R}^2$, and its objective value is normalised during surface construction.
We denote the resulting normalised objective value by $\widetilde{J}_i$.
Accordingly, we define a surface function $f:\mathbb{R}^2 \to \mathbb{R}$ through these samples by
\[
f(\boldsymbol y_i) = \widetilde{J}_i, \qquad i = 1,\ldots,N,
\]
that is, $f$ is known only through its sampled values on $\{\boldsymbol y_i\}$.
In implementation, a continuous visualisation surface is obtained from the scattered samples
$\{(\boldsymbol y_i, \widetilde{J}_i)\}$ using scattered-data interpolation; see \cite{scatteredInterpolant}.

\begin{figure}[!htp]
\centering
\includegraphics[width=0.8\textwidth, keepaspectratio, bb=192 168 3952 1226]{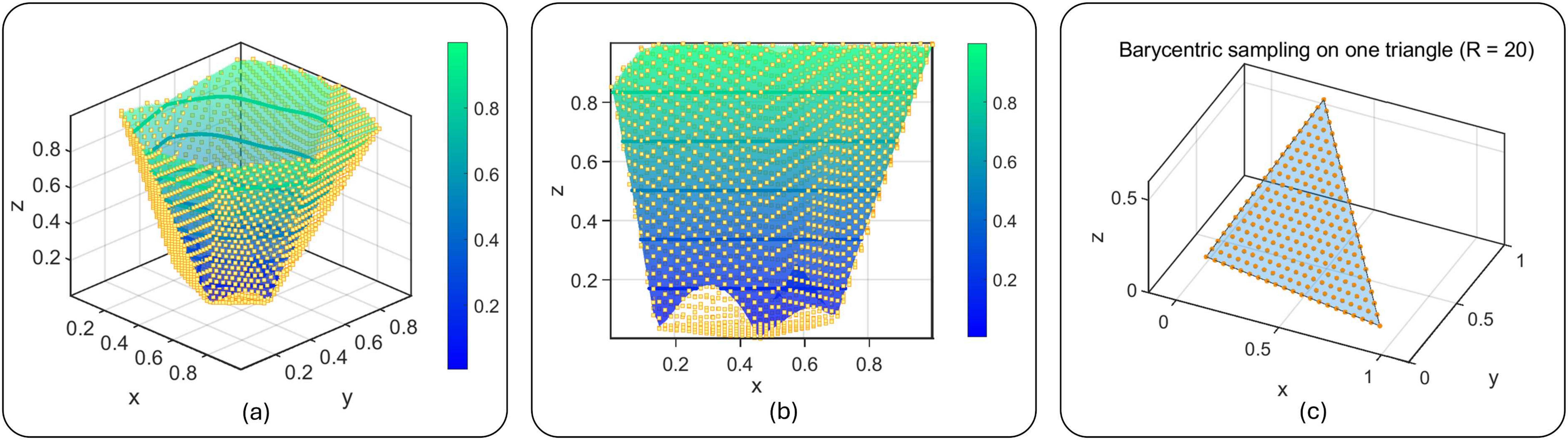}
\caption{Visualization for the compliance objective.
The blue surface is the original objective surface over $(x,y)$, and the gold mesh is the convex envelope obtained by the lower hull construction (Algorithm~\ref{alg:lower-hull-envelope}).
Panels (a) and (b) show two views of the same result.}
\label{exp:resultA1}
\end{figure}

To approximate the convex envelope of $f$, we consider the point set
$P = \{(\boldsymbol y_i, \widetilde{J}_i)\}_{i=1}^N \subset \mathbb{R}^3$
and compute its three-dimensional convex hull, from which the lower convex hull faces are extracted.
The convex hull construction follows standard algorithms; see \cite{convhulln}.
Let $F$ denote the set of triangular faces returned by the convex hull computation, and let
$F_{\downarrow} \subset F$ be the subset of lower convex hull faces selected according to
downward face normals (Algorithm~\ref{alg:lower-hull-envelope}).
On each retained face, barycentric samples of resolution $R$ are placed to obtain a set of
query points $\mathcal{Q} \subset \mathbb{R}^2$.
For each $\boldsymbol y \in \mathcal{Q}$, the value
$\widehat f_{\mathrm{ce}}(\boldsymbol y)$
is given by the $z$-coordinate of the corresponding barycentric sample point on $F_{\downarrow}$.
We then evaluate the pointwise distance between the visualisation surface and its lower convex hull
approximation at these barycentric samples.
Algorithm~\ref{alg:lower-hull-envelope} summarises the overall procedure.
Fig.~\ref{exp:resultA1}(a, b) show the original surface in blue together with gold barycentric
sample points that approximate its convex envelope, and Fig.~\ref{exp:resultA1}(c) zooms in on one
triangular face to highlight the barycentric sampling pattern.

\begin{algorithm}[!htbp]
\caption{Convex envelope approximation via lower convex hull}
\label{alg:lower-hull-envelope}
\begin{algorithmic}[1]
\Require Sample set $P=\{(\boldsymbol y_i,\widetilde{J}_i)\}_{i=1}^{N_{\mathrm{sample}}}\subset\mathbb{R}^3$; barycentric resolution $R$
\Ensure $\mathcal{Q}\subset\mathbb{R}^2$ and $\widehat{\mathcal{F}}=\{\widehat f_{\mathrm{ce}}(\boldsymbol y)\}_{\boldsymbol y\in\mathcal{Q}}$

\State $(F,V) \gets \texttt{convhulln}(P)$
\Comment{$V=\{\boldsymbol v_\ell\}\subset\mathbb{R}^3$ are convex hull vertices, $F$ are triangular faces (triples of vertex indices), see \cite{convhulln}};

\Statex \textit{Part I: select lower convex hull faces}
\State $F_{\downarrow}\gets \varnothing$
\For{each face $f=(i_1,i_2,i_3)\in F$}
  \State $p_1 \gets V[i_1]$, $p_2 \gets V[i_2]$, $p_3 \gets V[i_3]$
  \State $n\gets (p_2-p_1)\times(p_3-p_1)$
  \If{$n_z<0$}
    \State $F_{\downarrow}\gets F_{\downarrow}\cup\{f\}$
  \EndIf
\EndFor

\Statex \textit{Part II: barycentric sampling on the lower convex hull}
\State $\mathcal{Q}\gets \varnothing$, \; $\widehat{\mathcal{F}}\gets \varnothing$
\For{each face $f=(i_1,i_2,i_3)\in F_{\downarrow}$}
  \State $p_1 \gets V[i_1]$, $p_2 \gets V[i_2]$, $p_3 \gets V[i_3]$
  \For{$i=0$ to $R$}
    \For{$j=0$ to $R-i$}
      \State $a\gets i/R$, $b\gets j/R$, $c\gets 1-a-b$
      \State $p\gets c\,p_1+a\,p_2+b\,p_3$ \Comment{$p=(\boldsymbol y,\widehat f_{\mathrm{ce}}(\boldsymbol y))$}
      \State append $\boldsymbol y \gets (p_1,p_2)$ to $\mathcal{Q}$
      \State append $\widehat f_{\mathrm{ce}}(\boldsymbol y)\gets p_3$ to $\widehat{\mathcal{F}}$
    \EndFor
  \EndFor
\EndFor

\State \Return $\mathcal{Q}$ and $\widehat{\mathcal{F}}$
\end{algorithmic}
\end{algorithm}

\paragraph{Index construction}
We measure non-linearity by the vertical gap between the visualisation surface $f$
and its convex envelope approximation $\widehat f_{\mathrm{ce}}$, evaluated on the
barycentric sample set $\mathcal{Q}$.
Since the visualisation surface is constructed from normalised objective values,
both $f$ and $\widehat f_{\mathrm{ce}}$ take values in $[0,1]$ on $\mathcal{Q}$.
We therefore define the pointwise gap
\[
\Delta(\boldsymbol y)
= \big|\, f(\boldsymbol y)-\widehat f_{\mathrm{ce}}(\boldsymbol y)\,\big|
\ge 0, \qquad \boldsymbol y \in \mathcal{Q}.
\]
The scalar non-linearity index is defined as the mean gap,
\begin{equation}
\label{eq:index-mean-over-max}
\mathcal{I}_{\mathrm{NL}}
\;=\;
\frac{1}{|\mathcal{Q}|}\sum_{\boldsymbol y\in \mathcal{Q}} \Delta(\boldsymbol y)
\;\in\; [0,1].
\end{equation}
Since the vertical axis has already been normalised, the index is dimensionless.
Moreover, each retained lower-hull face is sampled with the same barycentric resolution,
so no additional area weighting is applied.

\section{Experiments}
\label{section3}

We systematically evaluate the proposed method through numerical experiments across multiple physics fields. In structural optimization, we consider compliance and stress problems, and additionally include the three field FPTO-based stress problem \cite{3field2024}. Performance is evaluated on 2D and 3D design domain. For heat transfer, the objectives include thermal compliance, maximum temperature, and temperature variance. For fluid flow, we test under varying Reynolds numbers. Finally, we conduct Four ablation studies to validate key design choices in our pipeline. First, we examine the dimensionality reduction step and show that our cosine similarity MDS embedding is well suited to distinguishing problem complexity.
Second, we study the effect of a sample group generated without fixing the objective gradient,
 whose best sample is treated as the reference solution, and compare the resulting visualization surfaces and the index $\mathcal{I}_{\mathrm{NL}}$. Third, we ablate the \emph{move limit} parameter $\eta_{\max}$ used in the fixed-gradient update (see Sec.~\ref{subsec:fixed_gradient}), and show how reducing the move limit mitigates excessive clustering or overlap of sampled points in the embedding. Fourth, we ablate the fixed-gradient strategy to verify its effectiveness in producing meaningful and stable sampling.
 All experiments are based on standard benchmark cases to ensure reproducibility and comparability. Each subsection presents the problem formulation, case setup, followed by visual and statistical results, with qualitative and quantitative analyses using ($\mathcal{I}_{\mathrm{NL}}$). To avoid distorting the cosine distance relationships among the initial points and to maintain consistency, the unfiltered design variables are used to construct the embedding space. In contrast, in the sampling process of Section~\ref{subsec:fixed_gradient} and in the presentation of standard optimization results, the material density filter is applied to the design variables in order to ensure physically meaningful material distributions, as shown in Figs.~\ref{setup}, \ref{distinct_local_optima}, \ref{problemsetting2}, \ref{problemsetting3}, and \ref{problemsetting4}.
 In addition, a small number of samples can yield abnormally large objective values that exceed the typical range. Such outliers can severely bias both the visualization surface and the resulting non-linearity evaluation. To mitigate this effect, we clip objective values using an objective-dependent upper bound: $10^4$ for compliance, $5\times 10^1$ for stress, $10^6$ for thermal compliance, $10^4$ for maximum temperature, $10^6$ for temperature variance, and $10^2$ for fluid energy dissipation. These bounds are chosen sufficiently high so that only anomalous values are affected. As a result, outliers are prevented from dominating the scale and obscuring the non-linearity in the key regions of interest. The numerical implementation of our TO frameworks is informed by established open-source codes and employs the original MMA code by Svanberg~\cite{MMA1987}; for detailed information regarding these computational aspects, please refer to~\ref{app1}.

\subsection{Structural optimization}
\label{subsec:struct_opt}

\paragraph{Common settings}
In the structural optimization experiments we consider three tasks, i.e., compliance minimization, maximum stress minimization, and maximum stress minimization via FPTO method. The first two are given in Eqs.~\eqref{eq:struct_opt_comp}-\eqref{eq:struct_opt_stress}. For FPTO method, interested readers can refer to \cite{3field2024} for full definitions, heuristic parameter choices, and extensions including body forces.

The compliance minimization task is described by
\begin{equation}
\label{eq:struct_opt_comp}
\begin{aligned}
\min_{\boldsymbol{\rho}\in[\rho_{\min},1]^n,\ \mathbf{u}}\quad
& J_{\mathrm{comp}}(\mathbf{u})=\mathbf{f}^{\top}\mathbf{u} \\
\text{s.t.}\quad
& \mathbf{K}(\boldsymbol{\rho})\,\mathbf{u}=\mathbf{f},\\
& \frac{\sum_{e=1}^{n} v_e\,\rho_e}{\sum_{e=1}^{n} v_e} \le \bar v .
\end{aligned}
\end{equation}
Where $\boldsymbol{\rho}=[\rho_1,\dots,\rho_n]^\top\!\in[\rho_{\min},1]^n$ are element density design variables with a small lower bound $\rho_{\min}>0$; $\mathbf{u}\in\mathbb{R}^{n_{\mathrm{dof}}}$ is the global displacement and $\mathbf{f}\in\mathbb{R}^{n_{\mathrm{dof}}}$ the global load; $\mathbf{K}(\boldsymbol{\rho})\in\mathbb{R}^{n_{\mathrm{dof}}\times n_{\mathrm{dof}}}$ is the assembled global stiffness matrix satisfying linear elasticity $\mathbf{K}(\boldsymbol{\rho})\mathbf{u}=\mathbf{f}$; $v_e$ is the volume of element $e$ and $\bar v\in(0,1]$ is the prescribed volume fraction; $J_{\mathrm{comp}}(\mathbf{u})=\mathbf{f}^{\top}\mathbf{u}$ denotes compliance (smaller values indicate stiffer designs).

\begin{figure}[!htp]
\centering
\includegraphics[width=0.8\textwidth, keepaspectratio, bb=137 111 709 405]{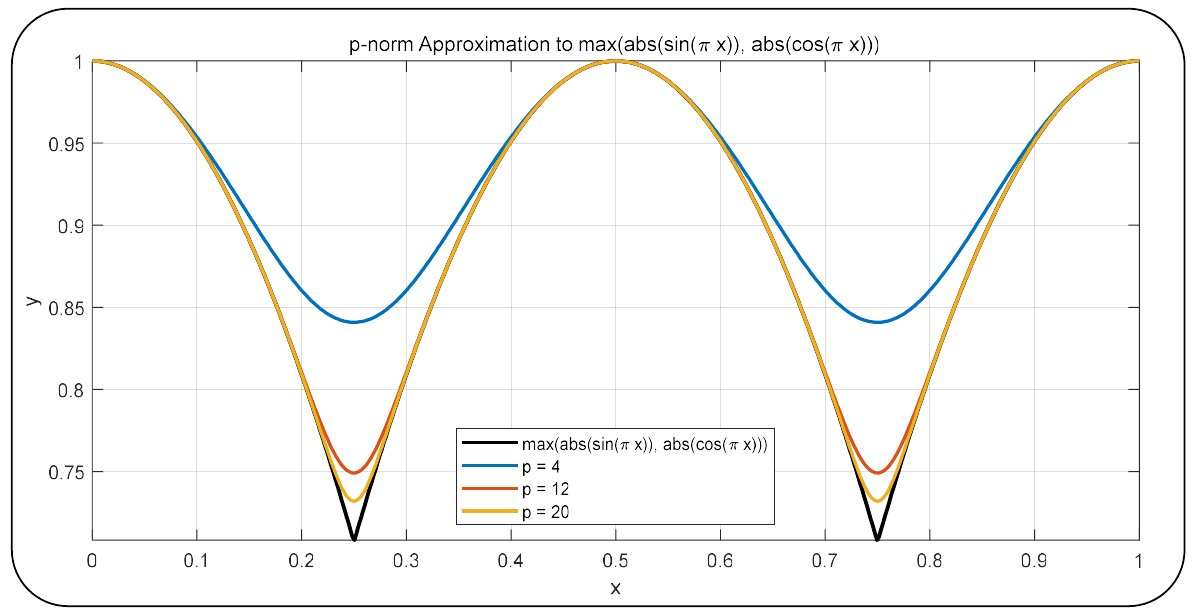}
\caption{\(p\)-norm approximation to a non-linear function.}
\label{pnormcompare}
\end{figure}

For maximum stress minimization we approximate the max operator with the $p$-norm aggregator:
\begin{equation}
\label{eq:struct_opt_stress}
\begin{aligned}
\min_{\boldsymbol{\rho}\in[\rho_{\min},1]^n,\ \mathbf{u}}\quad
& J_{\mathrm{stress}}(\mathbf{u},\boldsymbol{\rho})
  =\left(\sum_{e=1}^{n}  \,\sigma^{vm}_{e}(\mathbf{u},\boldsymbol{\rho})^{p}\right)^{\!1/p} \\[2pt]
\text{s.t.}\quad
& \mathbf{K}(\boldsymbol{\rho})\,\mathbf{u}=\mathbf{f},\\
& \frac{\sum_{e=1}^{n} v_e\,\rho_e}{\sum_{e=1}^{n} v_e} \le \bar v .
\end{aligned}
\end{equation}
where $\sigma^{vm}_{e}(\mathbf{u},\boldsymbol{\rho})$ is the von~Mises stress in element $e$, $p > 1$ controls the sharpness of the approximation and $J_{\mathrm{stress}}\to \max_e \sigma^{vm}_e$ as $p\to\infty$; the remaining symbols are as in Eq.~\eqref{eq:struct_opt_comp}. See Fig.~\ref{pnormcompare} for a 1D illustration: with $a(x)=|\sin(\pi x)|$ and $b(x)=|\cos(\pi x)|$ on $[0,1]$, the black curve plots $\max\{a(x),b(x)\}$ while the colored curves show the $p$-norm approximation $(a(x)^p+b(x)^p)^{1/p}$ for $p\in\{4,12,20\}$; larger $p$ tracks the maximum more closely but also increases non-linearity.

\begin{figure}[!htp]
\centering
\includegraphics[width=0.8\textwidth, keepaspectratio, bb=183 392 3897 1748]{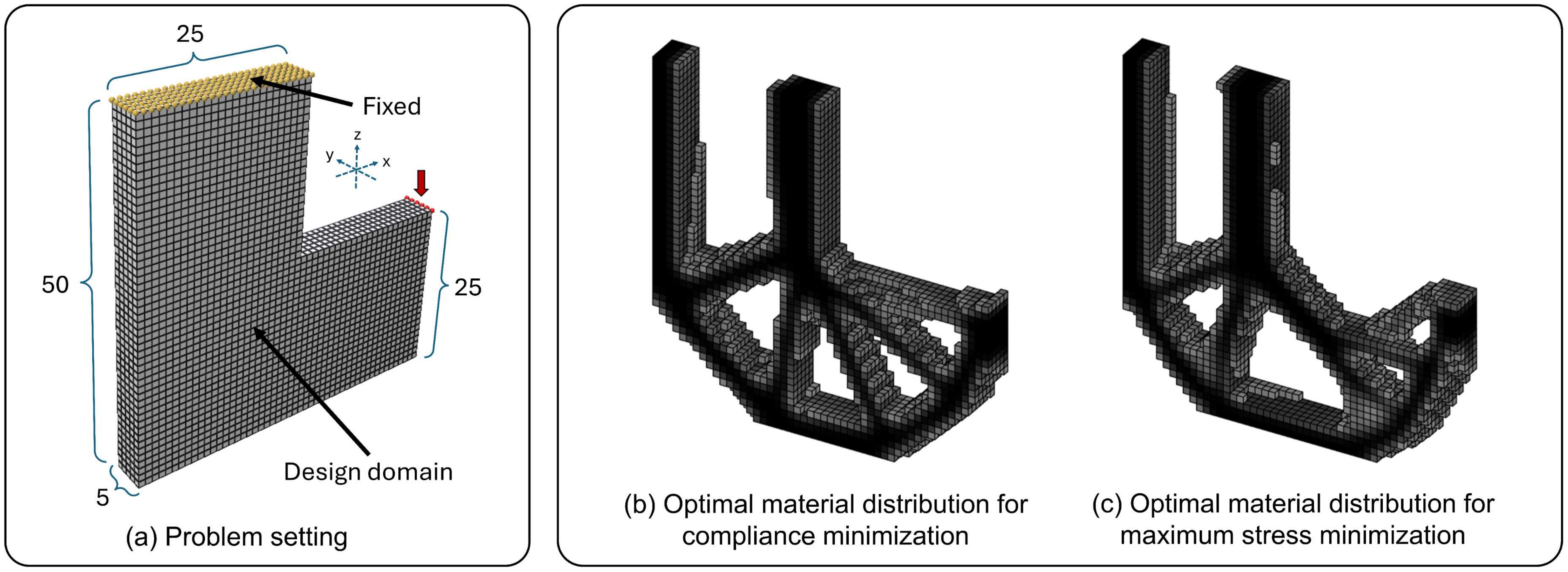}
\caption{Problem setting and optimization results.}

\label{problemsetting2}
\end{figure}

Fig.~\ref{setup} and Fig.~\ref{problemsetting2} illustrate the first and second problem settings used in our experiments, respectively. Fig.~\ref{problemsetting2}(a) shows the 3D design domain and boundary conditions, i.e., the dark gray voxels define the design domain, the golden vertices at the top end are fixed, and a downward unit load of \(1\,\text{N}\) is applied at the red vertices. 
Fig.~\ref{problemsetting2}(b) and (c) display representative optimal material distributions under the two objectives, i.e., (b) compliance minimization and (c) maximum stress minimization.

For the optimal solutions shown in Fig.~\ref{problemsetting2}(b, c), we use $\eta_{\max} = 0.2$. In the comparative experiments, we reduce the  $\eta_{\max}$ to 
0.01 after the fifth iteration. In the subsequent experiments, we start from different initial designs and record one sample every 20 iterations. This schedule ensures that, in every experiment, each group yields five samples before convergence and avoids excessive clustering of samples.

\begin{figure}[!htp]
\centering
\includegraphics[width=0.8\textwidth, keepaspectratio, bb=376 165 3885 2413]{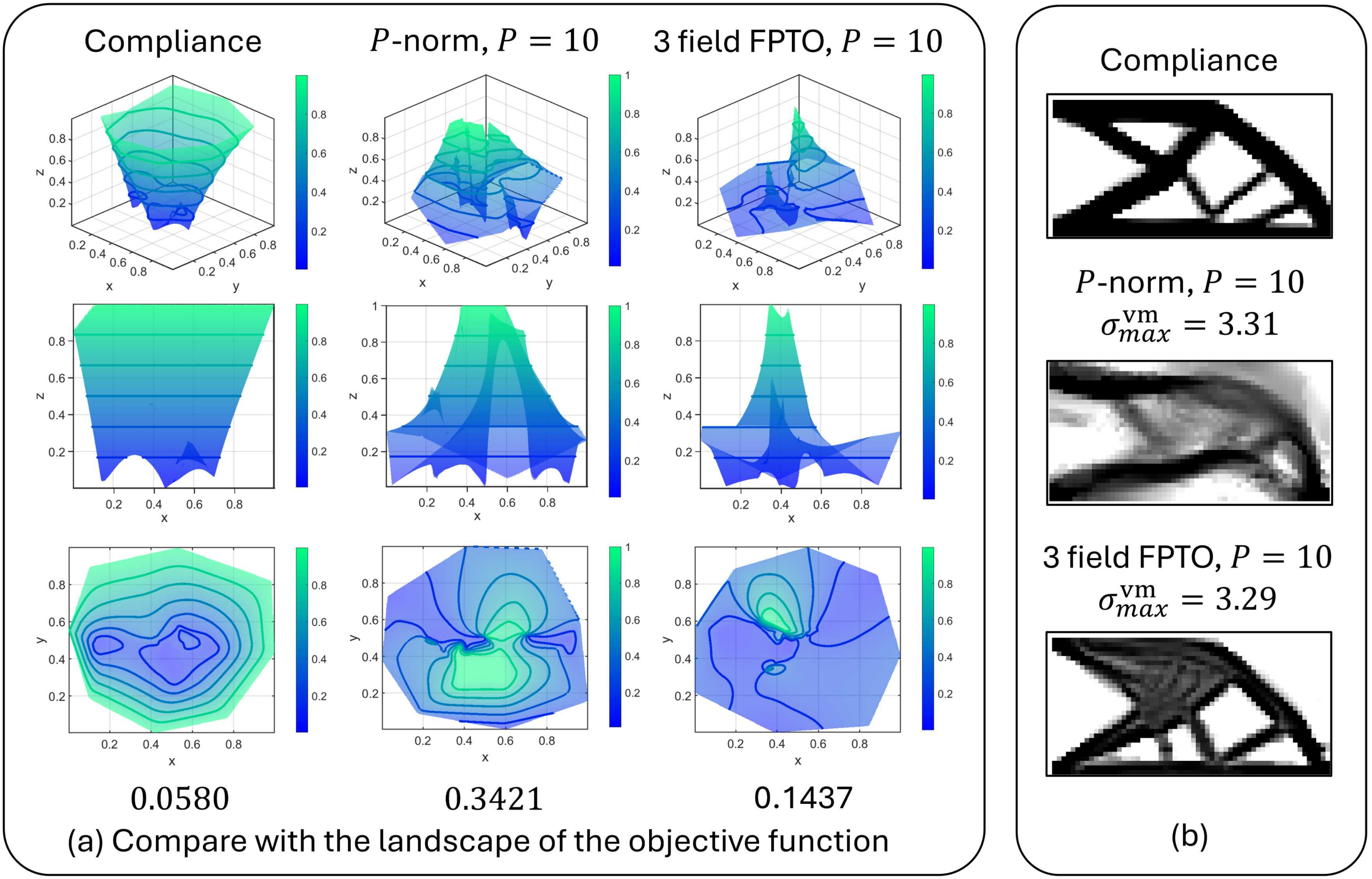}
\caption{(a) Landscape visualization for three objectives: compliance, $p$-norm stress ($p=10$), and FPTO method. Each column shows the objective surface and its contour maps on two-dimensional slices around 9 initial points; the number under each column is the non-linearity index $\mathcal{I}_{\mathrm{NL}}$. Colors encode the normalized objective value (higher is warmer). (b) Reference solutions computed without fixing the gradient.
}

\label{compare1}
\end{figure}

\paragraph{Comparison} 
Fig.~\ref{compare1}(a) shows the experiment under three different tasks. The left side visualizes the landscape of the objective function for the three tasks, i.e., compliance, $p$-norm stress with $p=10$, and FPTO with $p=10$, while Fig.~\ref{compare1}(b) shows the reference solutions obtained under identical algorithmic settings without fixing the gradient.

For compliance minimization the visualization surface is bowl shaped, with few peaks and valleys and a very smooth appearance, and sampling from the initial points along the fixed gradient directions only rarely increases the objective, which together indicates low non-linearity. The non-linearity index is $\mathcal{I}_{\mathrm{NL}}=0.0580$. For maximum stress minimization the landscape is peak dominated. During sampling the objective oscillates up and down, and the number of local peaks and valleys increases markedly. This behavior indicates high non-linearity, with the index rising to $\mathcal{I}_{\mathrm{NL}}=0.3421$. With FPTO, both the number and the amplitude (height/depth) of peaks and valleys decrease, indicating substantially fewer local optima and markedly reduced non-linearity. This reduction is consistent with \cite{3field2024}. By introducing a body force term and heuristically tuning the projection parameters, FPTO smooths the objective, which lowers non-linearity, improves convergence, and reduces sensitivity to parameters, leading to a more stable optimization.

Both optimal solutions use the same numerical settings as in Fig.~\ref{visualization}: optimizer, MMA; maximum allowable volume fraction, \(0.5\); maximum number of iterations, \(100\); and density filter radius, \(r = 2.5\) (in element size).

Correspondingly, the index decreases to $\mathcal{I}_{\mathrm{NL}}=0.1437$. As shown in Fig.~\ref{compare1}(b), under the same optimization settings, the compliance case has essentially converged, yielding a clear black-white layout with nearly no intermediate densities. By contrast, the maximum stress minimization case retains many gray elements, indicating incomplete convergence and strong non-linearity. After applying the FPTO method, the reference solution not only exhibits a clearer material layout, but also achieves a lower maximum von Mises stress ($\sigma_{\max}^{\mathrm{vm}}=3.29$) compared to the standard formulation ($\sigma_{\max}^{\mathrm{vm}}=3.31$), indicating an improvement in both solution performance and numerical robustness.

\begin{figure}[!htp]
\centering
\includegraphics[width=0.8\textwidth, keepaspectratio, bb=104 116 4090 2755]{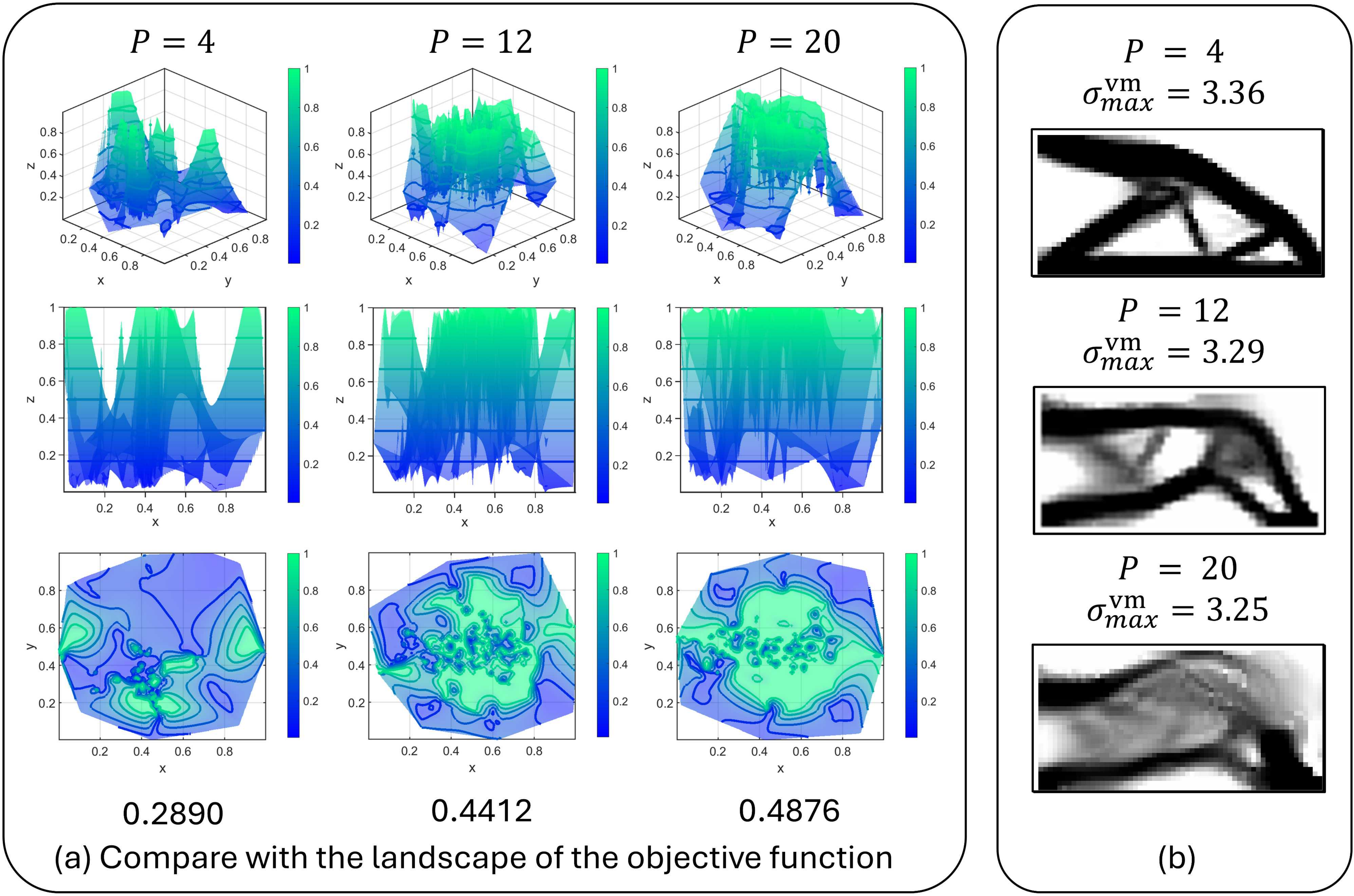}
\caption{(a) Landscape visualization for the $p$-norm stress objective with three different $p$ values ($p\in\{4,12,20\}$). Each column shows the visualization surface and its contour maps; the number under each column is the non-linearity index $\mathcal{I}_{\mathrm{NL}}$. Colors encode the normalized objective value (higher is warmer). (b) Reference solutions computed without fixing gradient.}

\label{compare2}
\end{figure}

Building on the objective-level comparison, we now examine the effect of the parameter $p$. To improve the accuracy of both visualization and quantification, we increase the number of initial points to 89. Fig.~\ref{compare2} reports results for $p\in\{4,12,20\}$. As $p$ grows, the landscape of the objective function becomes more rugged. Peaks increase in number and height. The non-linearity index rises from $\mathcal{I}_{\mathrm{NL}}=0.2890$ for $p=4$ to $\mathcal{I}_{\mathrm{NL}}=0.4412$ for $p=12$ and $\mathcal{I}_{\mathrm{NL}}=0.4876$ for $p=20$.

In the reference solutions, the case $p=4$ is nearly binary but has a higher max von Mises stress, $\sigma^{vm}_{\max}\approx 3.36$. The cases $p=12$ and $p=20$ retain some gray elements under the same settings, yet the max von Mises stress drops to about $3.29$ and $3.25$. These trends reflect that, as $p$ grows, the $p$-norm more closely approximates the $\max$ operator. As $p$ increases, non-linearity and optimization difficulty also increase. 
These observations are consistent with our quantitative and qualitative analyses and with the phenomenon in Fig.~\ref{pnormcompare}.

\begin{figure}[!htp]
\centering
\includegraphics[width=0.8\textwidth, keepaspectratio, bb=49 80 4024 2874]{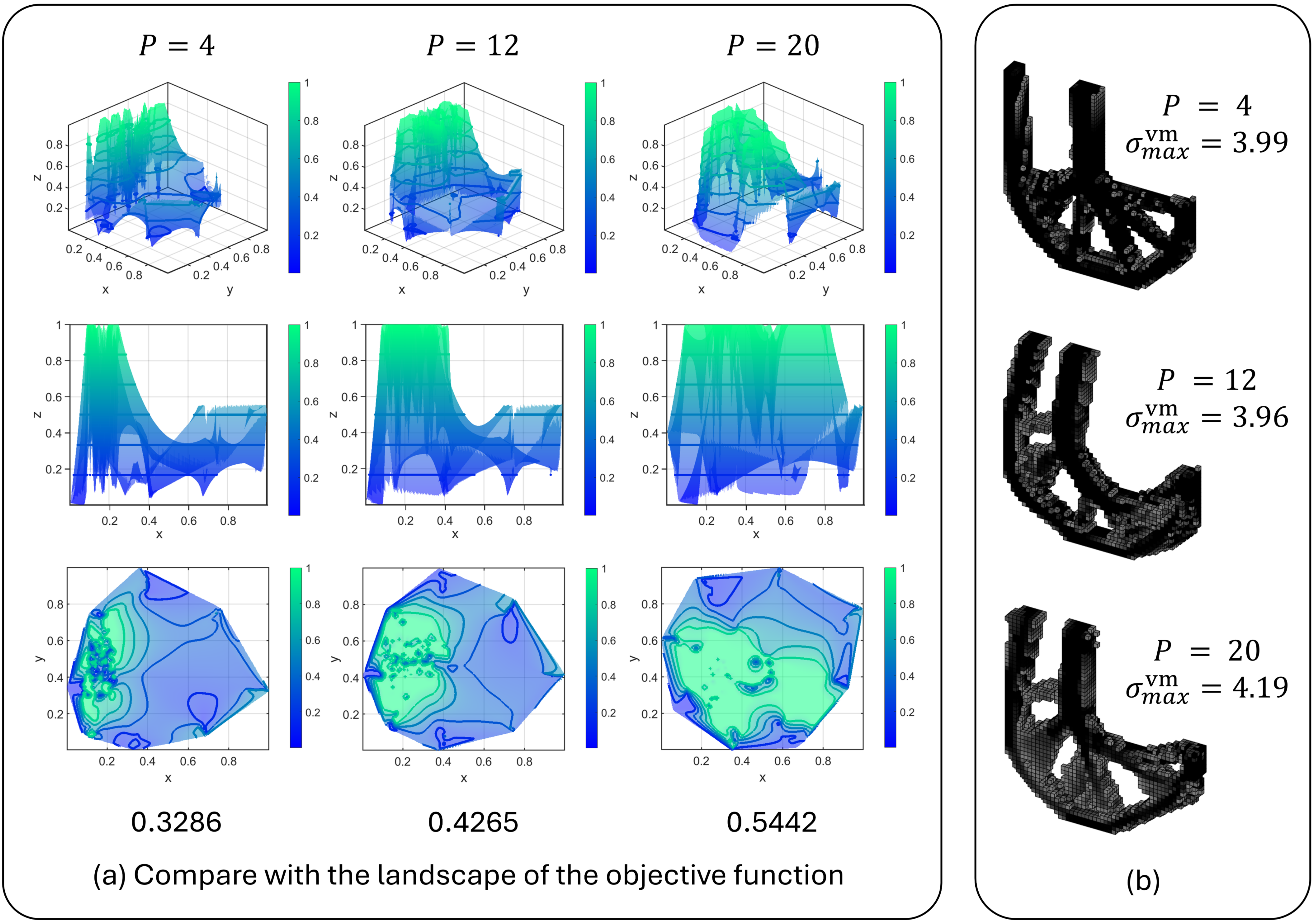}
\caption{(a) Landscape visualization for the $p$-norm stress objective with three different $p$ values ($p\in\{4,12,20\}$). Each column shows the visualization surface and its contour maps; the number under each column is the non-linearity index $\mathcal{I}_{\mathrm{NL}}$. Colors encode the normalized objective value (higher is warmer). (b) Reference solutions computed without fixing gradient.}

\label{compare2_3}
\end{figure}

To verify that our method remains applicable in a 3D design domain, we repeat the experiment in 3D with a thickness of 5 elements, using the same problem setting as in Fig.~\ref{problemsetting2}, and evaluate how different values of \(p\) affect the complexity of the problem.
Similar to the experiments in the 2D design domain, the non-linearity and optimization difficulty increase as \(p\) grows.
In particular, as shown in Fig.~\ref{compare2_3}, at \(p=12\) the reference solution performs much better than at \(p=20\), while the problem still exhibits substantially lower non-linearity.
This indicates that choosing \(p=12\) is more suitable for optimization than \(p=20\).

\subsection{Thermal optimization}
\label{subsec:thermal_opt}

\paragraph{Common settings}
We consider steady-state heat conduction with design $\boldsymbol{\rho}\in[\rho_{\min},1]^n$ and a volume fraction constraint $\sum_e v_e \rho_e / \sum_e v_e \le \bar v$.
The state equation is given by
$\mathbf{K}(\boldsymbol{\rho})\,\mathbf{T}=\mathbf{Q}$,
where $\mathbf{T}$ is the nodal temperature and $\mathbf{Q}$ collects the applied heat sources.

\medskip
\noindent\emph{Thermal compliance minimization.}
\begin{equation}
\label{eq:thermal_comp}
\begin{aligned}
\min_{\boldsymbol{\rho}\in[\rho_{\min},1]^n,\ \mathbf{T}}\quad
& J_{\mathrm{comp}}^{\mathrm{th}}(\mathbf{T})=\mathbf{Q}^{\top}\mathbf{T} \\
\text{s.t.}\quad
& \mathbf{K}(\boldsymbol{\rho})\,\mathbf{T}=\mathbf{Q}, \\
& \frac{\sum_{e=1}^{n} v_e\,\rho_e}{\sum_{e=1}^{n} v_e}\le \bar v .
\end{aligned}
\end{equation}
This is the thermal analogue of compliance in structural optimization.

\medskip
\noindent\emph{Maximum temperature minimization (via a \(p\)-norm).}
\begin{equation}
\label{eq:thermal_maxtemp}
\begin{aligned}
\min_{\boldsymbol{\rho}\in[\rho_{\min},1]^n,\ \mathbf{T}}\quad
& J_{\max T}(\mathbf{T})=\Big(\sum_{e=1}^{n} w_e\,T_e^{\,p}\Big)^{\!1/p} \\
\text{s.t.}\quad
& \mathbf{K}(\boldsymbol{\rho})\,\mathbf{T}=\mathbf{Q}, \\
& \frac{\sum_{e=1}^{n} v_e\,\rho_e}{\sum_{e=1}^{n} v_e}\le \bar v .
\end{aligned}
\end{equation}
where $T_e$ denotes an element-level temperature, $w_e>0$ are optional weights, usually equal to 1 in practice and $p>1$ controls the sharpness; as $p\to\infty$, $J_{\max T}\to \max_e T_e$.

\medskip
\noindent\emph{Temperature variance minimization.}
\begin{equation}
\label{eq:thermal_var}
\begin{aligned}
\min_{\boldsymbol{\rho}\in[\rho_{\min},1]^n,\ \mathbf{T}}\quad
& J_{\mathrm{var}}(\mathbf{T})
   = \sum_{e=1}^{n} w_e\,\big(T_e-\bar T\big)^2 \\[2pt]
\text{s.t.}\quad
& \mathbf{K}(\boldsymbol{\rho})\,\mathbf{T}=\mathbf{Q}, \\
& \frac{\sum_{e=1}^{n} v_e\,\rho_e}{\sum_{e=1}^{n} v_e}\le \bar v .
\end{aligned}
\end{equation}
where $\bar T=\big(\sum_e w_e T_e\big)/\big(\sum_e w_e\big)$ is the weighted mean temperature. 

\begin{figure}[!htp]
\centering
\includegraphics[width=0.8\textwidth, keepaspectratio, bb=170 288 3877 1416]{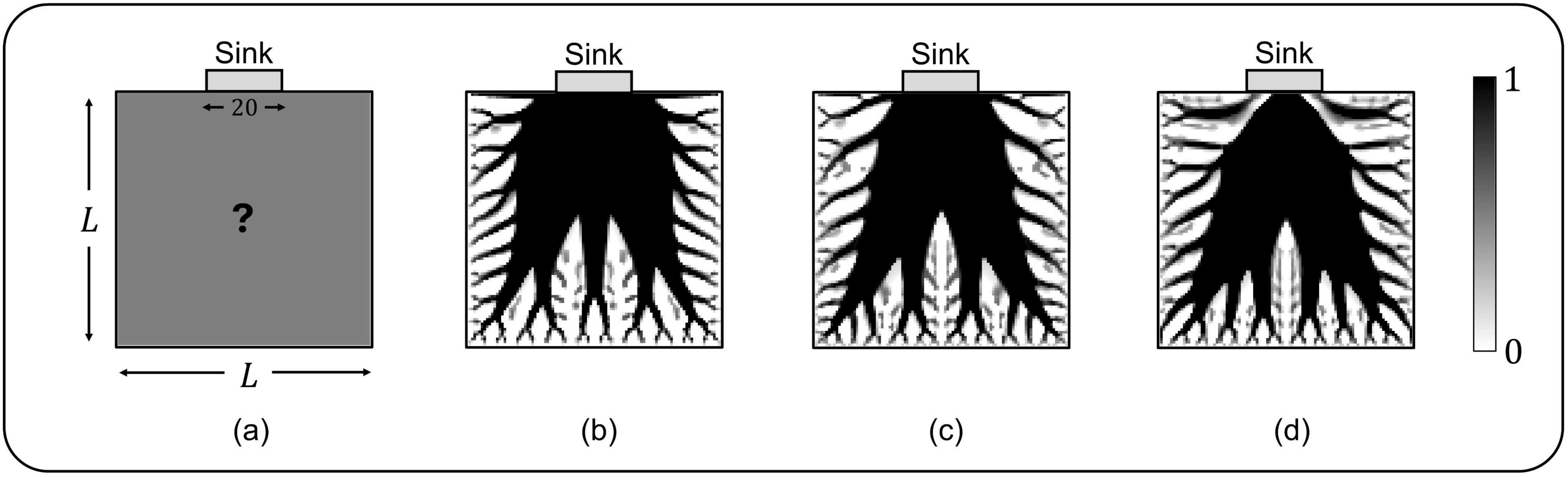}
\caption{Problem setting and optimization results for heat conduction. (a) Setup: the dark gray region is the design domain; the interface between the light gray “Sink” pad and the design domain at the top is fixed at ($T=0^\circ\mathrm{C}$). A uniform heat source ($q=-0.01$) is applied throughout the design domain. (b) Optimal solution for thermal compliance minimization. (c) Maximum temperature minimization. (d) Temperature variance minimization. The color bar from 0 to 1 denotes material density. All other settings are identical to those in Fig.~\ref{setup}.
}
\label{problemsetting3}
\end{figure}

The problem setting is shown in Fig.~\ref{problemsetting3}. The computational domain is a square with a mesh resolution of $L\times L$ cells. The top boundary contains a centered isothermal heat sink segment spanning 20 cells, modeled as a Dirichlet condition
$T=0^\circ\mathrm{C}$. The remainder of the boundary (including the rest of the top edge) is adiabatic. In the heat transfer problem, to demonstrate the robustness of the proposed method with respect to different parameter configurations, we adjust the maximum volume constraint to \(0.6\). All other parameters follow the settings above. For the optimal solutions shown in Fig.~\ref{problemsetting3}(b)-(d), we use a  $\eta_{\max}$ of 
0.2. In the comparative experiments, consistent with the previous setup, we reduce the  $\eta_{\max}$ to 
0.01 after the fifth iteration.

\begin{figure}[!htp]
\centering
\includegraphics[width=0.8\textwidth, keepaspectratio, bb= 239 105 4033 2632]{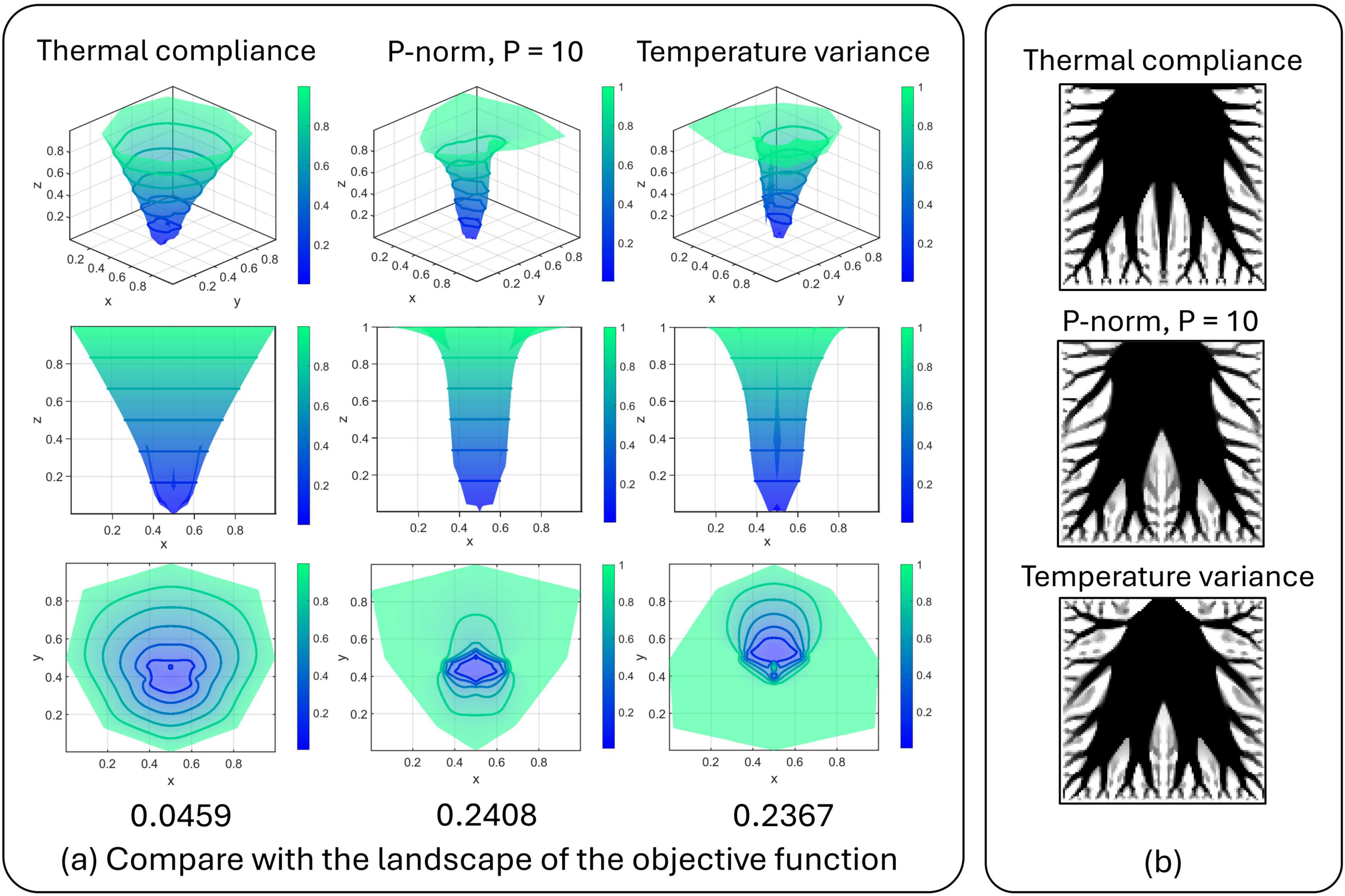}
\caption{(a) Visualization surfaces for three thermal objectives: thermal compliance, maximum temperature via a $p$-norm ($p=10$), and temperature variance. Each column shows the visualization surface and its contour maps; the value under each column is the non-linearity index $\mathcal{I}_{\mathrm{NL}}$. Colors encode the objective value (higher is warmer). (b) Reference solutions computed without fixing gradient.}
\label{compare3}
\end{figure}

\paragraph{Comparison}
Fig.~\ref{compare3} summarises the results for the heat transfer problem.
Across the three thermal objectives, the visualisation surfaces all exhibit a clear valley-like structure, indicating that these problems are generally less non-linear than the structural TO cases.
Among them, thermal compliance case yields a surface that is closest to its convex envelope, reflected by the smallest non-linearity index $\mathcal{I}_{\mathrm{NL}}=0.0459$.
The other two objectives, namely the maximum temperature objective approximated by a $p$-norm ($p=10$) and the temperature variance objective, produce visibly more complex surfaces, with larger indices $\mathcal{I}_{\mathrm{NL}}=0.2408$ and $\mathcal{I}_{\mathrm{NL}}=0.2301$, respectively.

This trend is consistent with the reference solutions shown on the right.
For all three objectives, elements with intermediate densities are mainly concentrated near the boundary regions.
Compared with the latter two objectives, thermal compliance results in fewer such intermediate elements and a layout that is closer to a binary distribution.
These observations are in agreement with the quantitative non-linearity indices and confirm the comparatively low non-linearity of the heat transfer problems considered here.

\begin{figure}[!htp]
\centering
\includegraphics[width=0.75\textwidth, keepaspectratio, bb=109 97 4041 2841]{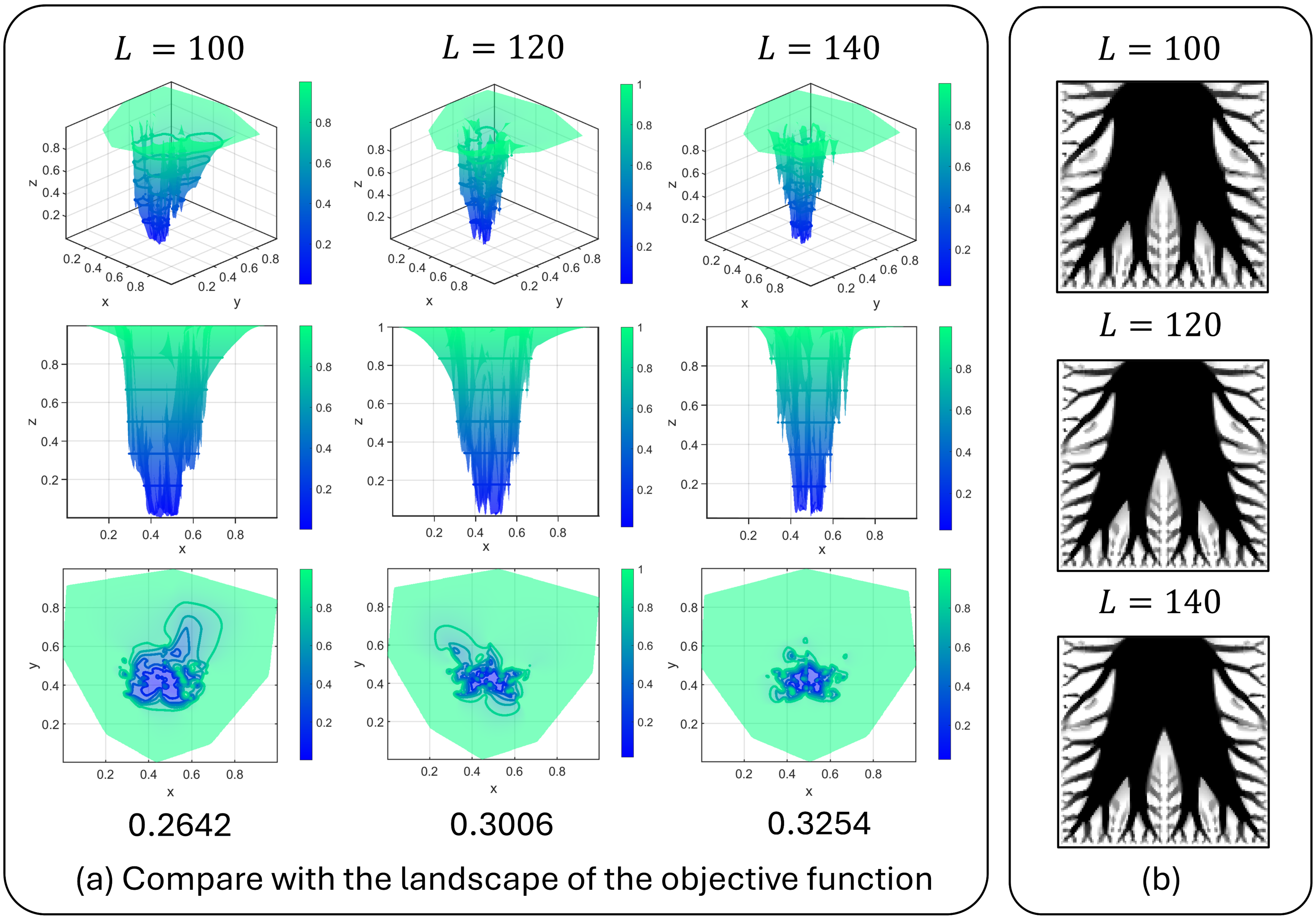}
\caption{(a) Visualization surfaces for the maximum temperature objective at $L=100, 120, 140$. Each column shows the visualization surface and its contour maps surrounded by 89 initial points; the value under each column is the non-linearity index $\mathcal{I}_{\mathrm{NL}}$. Colors encode the objective value (higher is warmer). (b) Reference solutions computed without fixing gradient.}

\label{compare4}
\end{figure}

We next examine the impact of the mesh resolution \(L\) (see Fig.~\ref{compare4}). 
For \(L\in\{100,120,140\}\), the non-linearity index increases from 
\(\mathcal{I}_{\mathrm{NL}}=0.2642\) (\(L=100\)) to 
\(\mathcal{I}_{\mathrm{NL}}=0.3006\) (\(L=120\)) and 
\(\mathcal{I}_{\mathrm{NL}}=0.3254\) (\(L=140\)). 
As the mesh is refined, the visualization surface develops deeper valleys and more intricate contour patterns, indicating a gradual increase in non-linearity and problem complexity.

On the reference solution side, the designs remain largely unchanged with varying mesh resolution $L$, indicating that the heat transfer problem is inherently less non-linear than structural TO problems.
Even in this low nonlinearity and mesh insensitive situation, our method yields a reasonable and consistent non-linearity measure, demonstrating its robustness.



\subsection{Fluid-flow optimization}
\label{subsec:flow_opt}

\begin{figure}[!htp]
\centering
\includegraphics[width=0.8\textwidth, keepaspectratio, bb=208 82 3802 955]{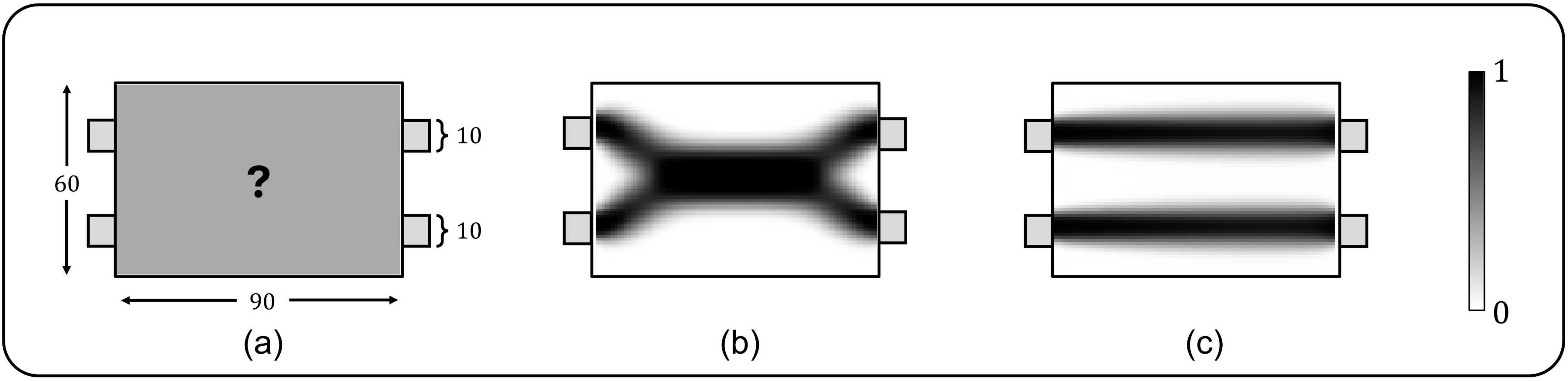}
\caption{Flow topology optimization results. (a) Problem setup: the dark gray region is the design domain. (b) Result at Reynolds number ($Re=0.17$). (c) Result at ($Re=133$). The color bar from $0$ to $1$ indicates material density.}

\label{problemsetting4}
\end{figure}

\paragraph{Common settings}
We consider steady, incompressible laminar flow on the whole domain using the Brinkman penalized Navier Stokes model. Let $\boldsymbol{\rho}\in[\rho_{\min},1]^n$ be the elementwise design field. The state variables are velocity $\mathbf{u}$ and pressure $p$. With a spatially varying Brinkman factor $\alpha(\boldsymbol{\rho})\ge 0$, the governing equations are
\begin{equation}
\label{eq:flow_state}
\begin{aligned}
& \rho_f(\mathbf{u}\!\cdot\!\nabla)\mathbf{u}
 - \mu\,\nabla\!\cdot\!\big(\nabla\mathbf{u}+\nabla\mathbf{u}^{\!\top}\big)
 + \nabla p
 + \alpha(\boldsymbol{\rho})\,\mathbf{u}
 = \mathbf{0}\quad\text{in }\Omega,\\
& \nabla\!\cdot\!\mathbf{u}=0\quad\text{in }\Omega,
\end{aligned}
\end{equation}
which extends the equations to both fluid and solid regions. Here, $\rho_f$ denotes the fluid density; and $\mu$ is the dynamic viscosity; see \cite{fieldTO} for details.

\medskip
\noindent\emph{Dissipated energy minimization.}
Following the classical setting, we minimize the total energy dissipation:
\begin{equation}
\label{eq:flow_obj}
\begin{aligned}
\min_{\boldsymbol{\rho}\in[\rho_{\min},1]^n,\ \mathbf{u},\,p}\quad
& J_{\mathrm{diss}}(\mathbf{u},\boldsymbol{\rho})
= \int_{\Omega}
\Big[
\tfrac{\mu}{2}\,(\nabla\mathbf{u}+\nabla\mathbf{u}^{\!\top}):(\nabla\mathbf{u}+\nabla\mathbf{u}^{\!\top})
+\alpha(\boldsymbol{\rho})\,\mathbf{u}\!\cdot\!\mathbf{u}
\Big]\; \mathrm{d}V \\[2pt]
\text{s.t.}\quad
& \text{\eqref{eq:flow_state}},\qquad
V \;\le\; \bar V ,
\end{aligned}
\end{equation}
with the relative fluid volume
\begin{equation}
\label{eq:flow_volume_cont}
V \;=\; \frac{1}{|\Omega|}\int_{\Omega}\gamma(\mathbf{x})\,\mathrm{d}V,
\qquad
\gamma(\mathbf{x})=
\begin{cases}
1, & \mathbf{x}\in\Omega_f,\\
0, & \mathbf{x}\in\Omega_s.
\end{cases}
\end{equation}
Here \(\Omega\subset\mathbb{R}^2\) denotes the entire computational domain, which is decomposed into the fluid region \(\Omega_f\) and the solid region \(\Omega_s\), with \(\Omega = \Omega_f \cup \Omega_s\) and \(\Omega_f \cap \Omega_s = \emptyset\). The characteristic function \(\gamma\) in \eqref{eq:flow_volume_cont} selects the fluid part of the domain.

\begin{figure}[!htp]
\centering
\includegraphics[width=0.8\textwidth, keepaspectratio, bb=55 52 4027 2841]{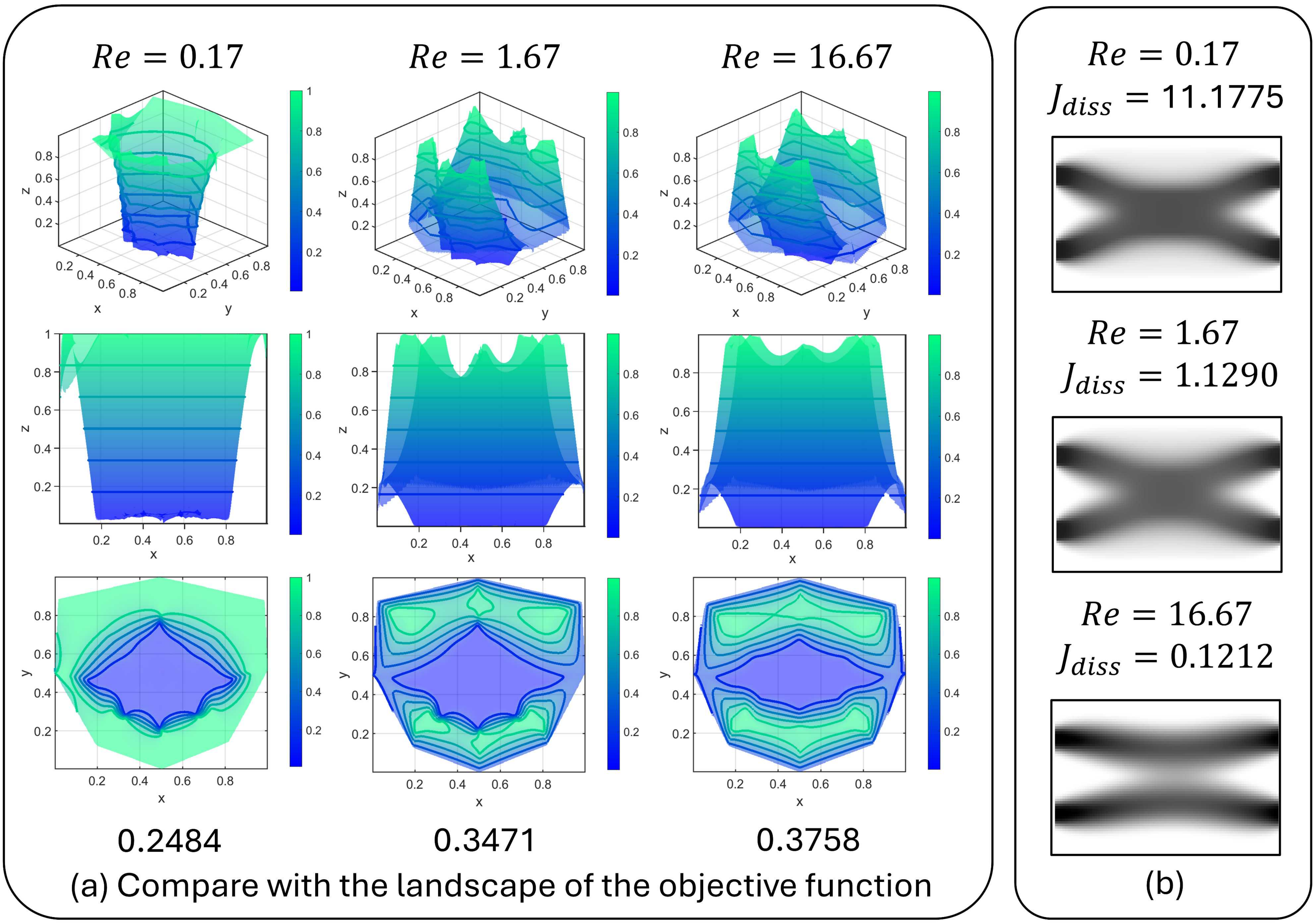}
\caption{(a) Visualization surfaces for the flow problem at three Reynolds numbers, $\mathrm{Re}\in\{0.17,\,1.67,\,16.67\}$. Each column shows the visualization surface and its contour maps surrounded by nine initial points; the value under each column is the non-linearity index $\mathcal{I}_{\mathrm{NL}}$. Colors encode the objective value (higher is warmer). (b) Reference solutions computed without fixing gradient.}

\label{compare5}
\end{figure}

Fig.~\ref{problemsetting4} shows the setup used in our experiments. The dark gray rectangle indicates the design part of \(\Omega\); the short light gray stubs on the left and right act as inlet and outlet channels, and all external walls are no slip. We solve the steady laminar problem \eqref{eq:flow_state} with the Brinkman extension and compare different Reynolds numbers \(Re\) by varying the characteristic inlet velocity. Representative optimized layouts are shown on the right.

\paragraph{Comparison}
Fig.~\ref{compare5} summarizes the results across three Reynolds numbers,
$Re \in \{0.17,\,1.67,\,16.67\}$. The geometry, boundary conditions,
and volume fraction constraint ($\bar v = 0.34$) follow
Fig.~\ref{problemsetting4}. From nine distinct initial points we perform
sampling to construct the visualization surfaces. In panel (a), each column
shows the surface, its side view, and the corresponding contour map; the
value printed below each column is the non-linearity index
$\mathcal{I}_{\mathrm{NL}}$.

As the Reynolds number increases, the visualization surfaces evolve from a predominantly valley-like structure to a more mountainous landscape with multiple peaks and ridges, indicating a gradual increase in problem non-linearity.
Consistently, the non-linearity index $\mathcal{I}_{\mathrm{NL}}$ increases with the Reynolds number, quantitatively capturing this progressive growth in complexity.

\subsection{Ablation studies}
\label{sec:ablationaaa}
We present a robust visualization framework together with a geometrically meaningful quantitative index. In this section, we conduct four ablation studies, where we replace one key step or parameter of the proposed pipeline while keeping all other steps unchanged, and then examine the resulting changes in the visualization surfaces and in $\mathcal{I}_{\mathrm{NL}}$. This ablation protocol isolates the contribution of each design choice and clarifies the role played by the corresponding component within the overall framework. We consider four ablations: (i) the dimensionality reduction method used for embedding, (ii) The effect of adding a group sample points generated by a standard optimization run without fixing the gradient, And (iii) the  $\eta_{\max}$ used during sampling after the gradient is fixed. (iv) whether the objective gradient is fixed during sampling.

\begin{figure}[!htp]
\centering
\includegraphics[width=0.8\textwidth, keepaspectratio, bb= 236 291 3835 3898]{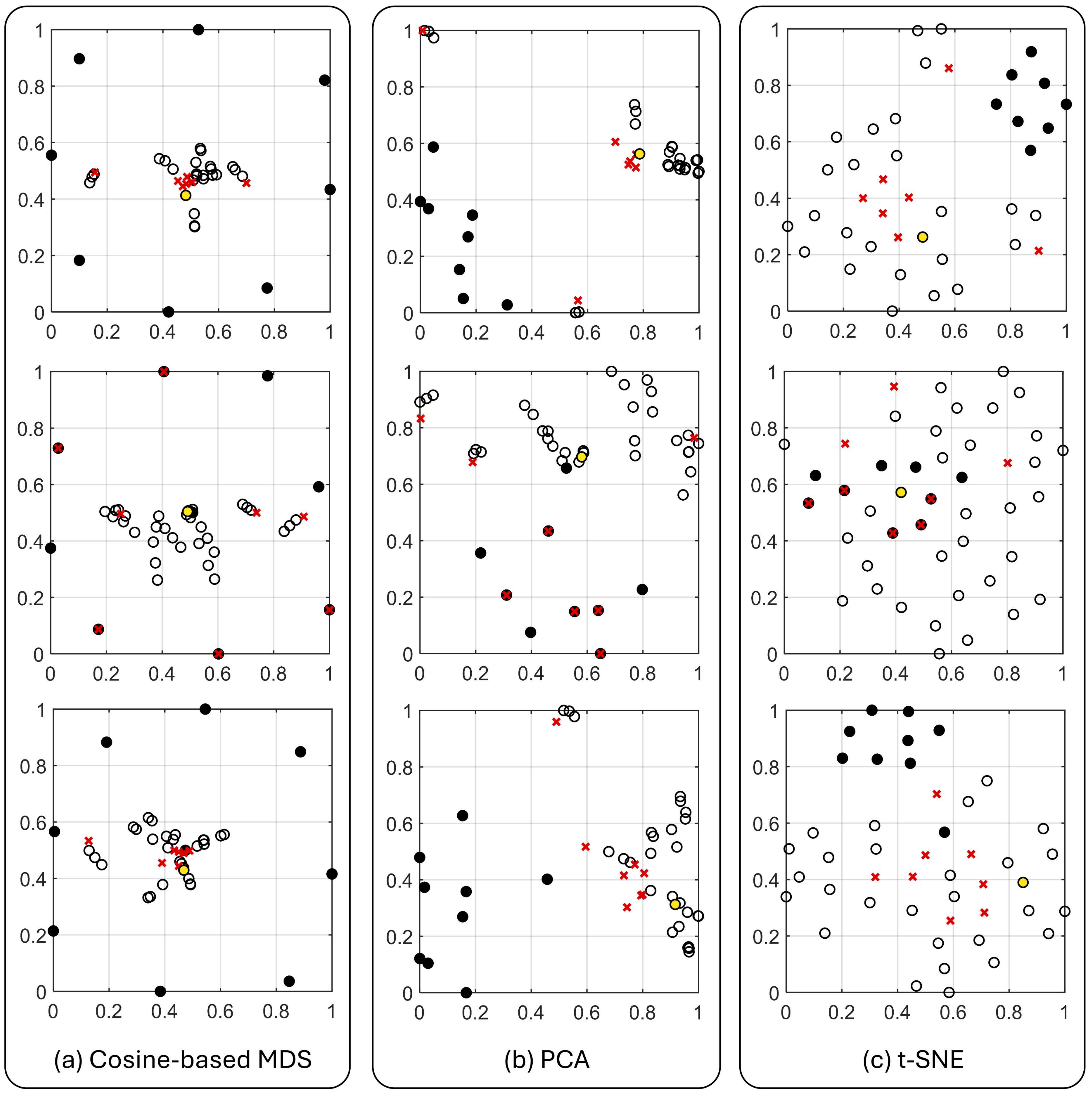}
\caption{Embeddings of the same sampling datasets with different dimensionality-reduction methods: (a) cosine-based MDS, (b) PCA, and (c) t-SNE. The rows correspond to different physical problems: the first row represents compliance minimization, the second row maximum stress minimization, and the third row the three field FPTO problem. All parameters are consistent with those in Fig.~\ref{compare1}. Black filled circles mark the starting points of each group, gold filled circles denote reference solutions that we treat as global optima, and red crosses indicate the per-group best sample.}
 \label{aaa1}
\end{figure}

\paragraph{Embedding method}
Fig.~\ref{aaa1} compares cosine-based MDS, PCA, and t-SNE \cite{tsne} using the same sampling datasets.
Across all embeddings, a consistent trend can be observed between the degree of clustering among sample points and the underlying non-linearity of the TO problem.
In particular, for maximum stress minimization problem, the embedded samples are noticeably more dispersed than those of compliance minimization, indicating a higher level of non-linearity.
After applying the three-field FPTO formulation, the samples exhibit a substantially increased level of clustering, suggesting that the non-linearity of the problem is significantly reduced.
However, this comparison also highlights a limitation, i.e., relying solely on the visual concentration or dispersion of embedded samples is insufficient for a rigorous and quantitative assessment of non-linearity.
This observation directly motivates our non-linearity index based on the gap between the visualization surface and its lower convex hull.
Since the objective function values are not involved in the dimensionality reduction process, incorporating them at the quantification stage avoids information loss and enables a more accurate and robust characterization of non-linear complexity.

\paragraph{Sampling strategy and reference solution}

To better demonstrate the proposed pipeline, we include five additional samples in all numerical experiments and visualizations. These samples are generated by a standard optimization run without fixing the gradient, and the best sample within this group is denoted as the reference solution.

\begin{figure}[!htp]
\centering
\includegraphics[width=0.75\textwidth, keepaspectratio, bb=303 200 3871 2671]{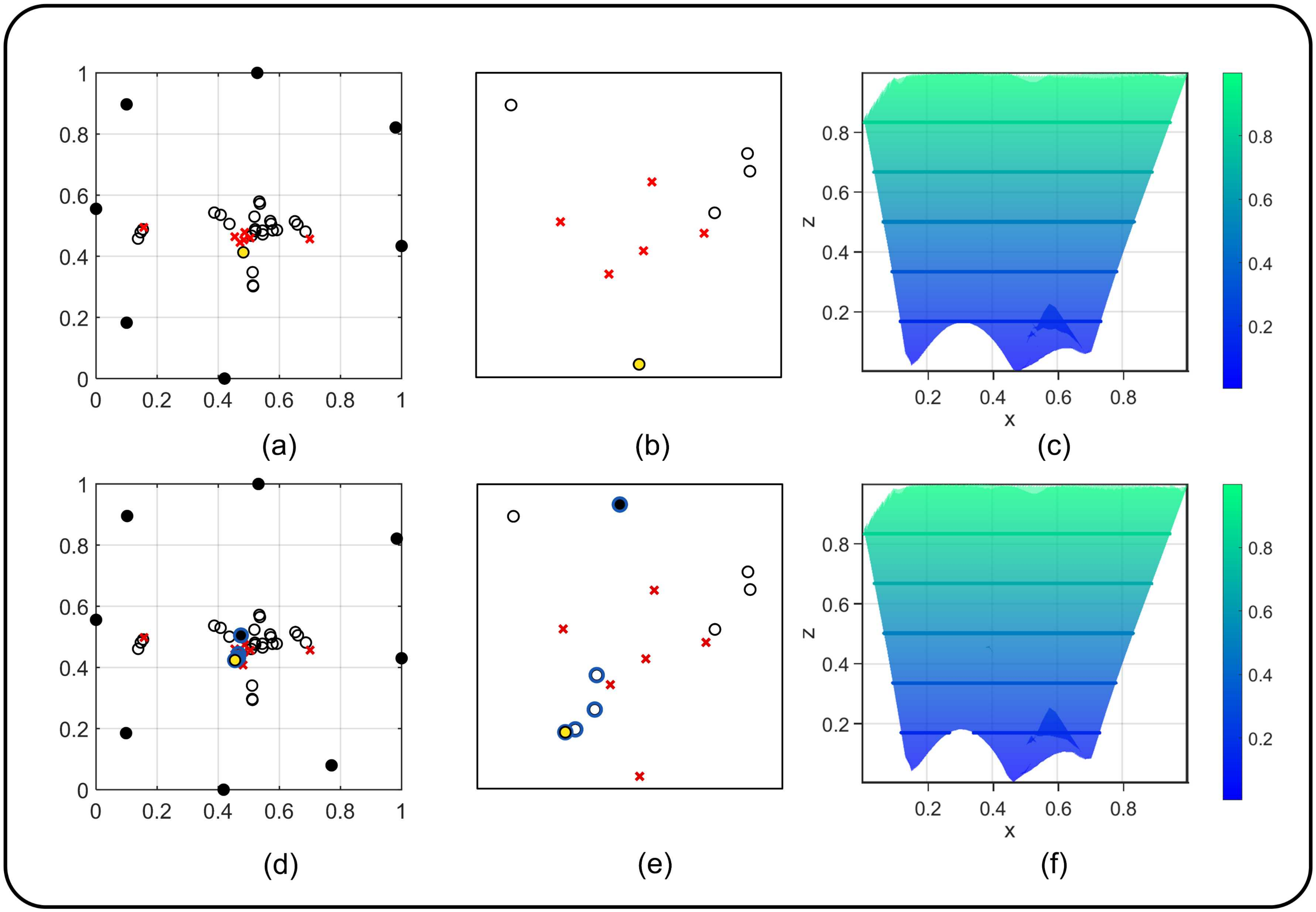}
\caption{Ablation on multi-start fixed-gradient sampling and the effect of adding a small sample group generated without fixing the gradient (whose best sample is treated as the reference solution). The corresponding numerical example is presented in the first column of Fig.~\ref{compare1}.}

\label{ablation2}
\end{figure}

Fig.~\ref{ablation2} evaluates the proposed multi-start fixed-gradient sampling strategy.
The top row employs eight starting points, all of which are sampling with fixed objective gradients.
The bottom row additionally introduces a ninth starting point located at the center of the embedding, indicated by a solid black marker, from which five additional samples are generated using a standard optimization run without fixing the gradient. These additional samples are highlighted by blue circular outlines.
The first column shows the embedding of all sampled designs, while the second column presents a local zoom-in of the corresponding region to better illustrate the distribution of samples. The third column depicts the reconstructed visualization surface. As can be observed, the resulting visualization surface remain almost unchanged after adding the additional sample group, demonstrating that the additional ninth sample group, introduced solely for demonstration purposes, does not lead to a noticeable change in the visualization surface.

Let $\mathcal{I}_{\mathrm{NL}}^{(i)}$ and $\widetilde{\mathcal{I}}_{\mathrm{NL}}^{(i)}$ denote the non-linearity index obtained in the $i$-th experiment without and with this additional sample group, respectively.
To quantitatively assess the influence of introducing this additional sample group, we define the absolute percentage change as
\begin{equation}
\Delta_i
=
\left|
\frac{\widetilde{\mathcal{I}}_{\mathrm{NL}}^{(i)}-\mathcal{I}_{\mathrm{NL}}^{(i)}}
{\mathcal{I}_{\mathrm{NL}}^{(i)}}
\right|\times 100\%.
\end{equation}
We further report the mean absolute percentage change over all $K$ experiments,
\begin{equation}
\overline{\Delta}
=
\frac{1}{K}\sum_{i=1}^{K}\Delta_i ,
\end{equation}
where $K$ is the total number of experiments considered.

As summarized in Table~\ref{tab:additional_group_ablation}, the inclusion of this additional sample group results in only marginal variations in $\mathcal{I}_{\mathrm{NL}}$ across all experiments.
The resulting mean absolute percentage change is $\overline{\Delta}=0.594\%$.

\begin{table}[!htbp]
\centering

\begin{tabular}{c|ccc}
\hline
Experiment & Case A & Case B & Case C \\
\hline
\multicolumn{4}{c}{\textbf{Without additional samples}} \\
\hline
1 & 0.0588 & 0.3574 & 0.1435 \\
2 & 0.2895 & 0.4415 & 0.4891 \\
3 & 0.3582 & 0.4776 & 0.5481 \\
4 & 0.0385 & 0.2434 & 0.2389 \\
5 & 0.2632 & 0.2964 & 0.3294 \\
6 & 0.2526 & 0.3477 & 0.3767 \\
\hline
\multicolumn{4}{c}{\textbf{With additional samples}} \\
\hline
1 & 0.0580 & 0.3421 & 0.1437 \\
2 & 0.2890 & 0.4412 & 0.4876 \\
3 & 0.3286 & 0.4265 & 0.5442 \\
4 & 0.0459 & 0.2408 & 0.2367 \\
5 & 0.2642 & 0.3006 & 0.3254 \\
6 & 0.2484 & 0.3471 & 0.3758 \\
\hline
\end{tabular}
\caption{Non-linearity index $\mathcal{I}_{\mathrm{NL}}$ with and without the additional sample group.}
\label{tab:additional_group_ablation}

\end{table}

\begin{figure}[!htp]
\centering
\includegraphics[width=0.8\textwidth, keepaspectratio, bb=69 236 3885 4037]{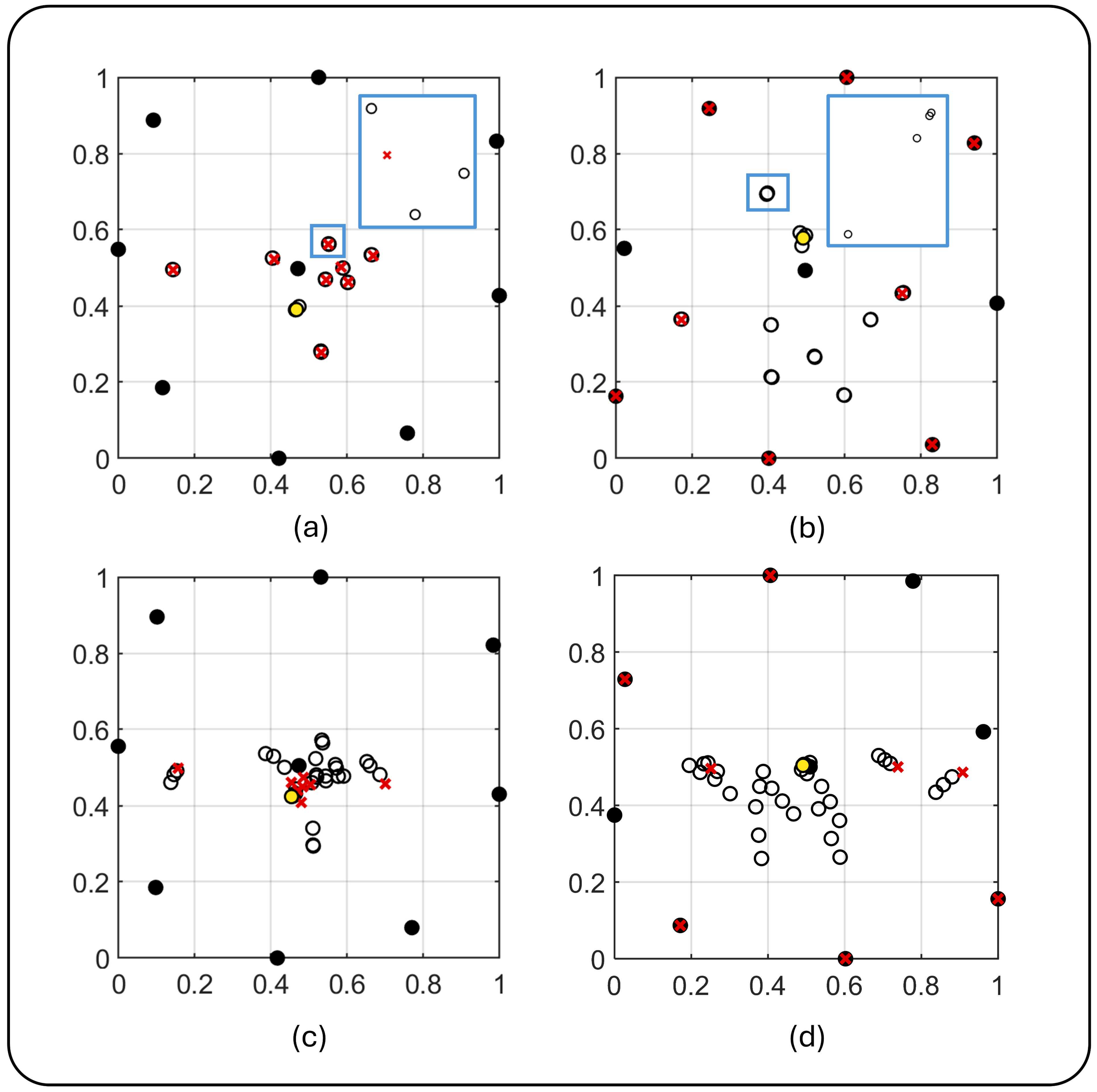}
\caption{Ablation study on the choice of $\eta_{\max}$ after fixing the objective gradient. Panels (a) and (c) show the 2D embeddings for the compliance minimization problem, while panels (b) and (d) correspond to the maximum stress minimization problem. The first row uses a fixed $\eta_{\max}=0.2$ throughout, whereas in the second row $\eta_{\max}$ is reduced to $0.01$ after fixing the gradient. Blue boxes indicate locally magnified regions. Without reducing $\eta_{\max}$, many sampled points become clustered or even overlap, which effectively reduces the number of distinct samples. In contrast, decreasing $\eta_{\max}$ leads to better-separated samples and improves the effective sampling quality. The corresponding numerical example is presented in the first and second column of Fig.~\ref{compare1}.}
\label{ablation3}
\end{figure}

\paragraph{Ablation on the parameter $\eta_{\max}$}
We further conduct a third ablation on the  $\eta_{\max}$ used after the gradient is fixed. We compare results under different  $\eta_{\max}$ and examine the resulting changes in the embedding space (Fig.~\ref{ablation3}). This ablation motivates our choice to reduce the  $\eta_{\max}$ after fixing the gradient, which helps prevent the sampled points from becoming overly clustered or even overlapping, and improves the comparability of intermediate samples.

\begin{figure}[!htp]
\centering
\includegraphics[width=0.8\textwidth, keepaspectratio, bb=95 236 3898 4037]{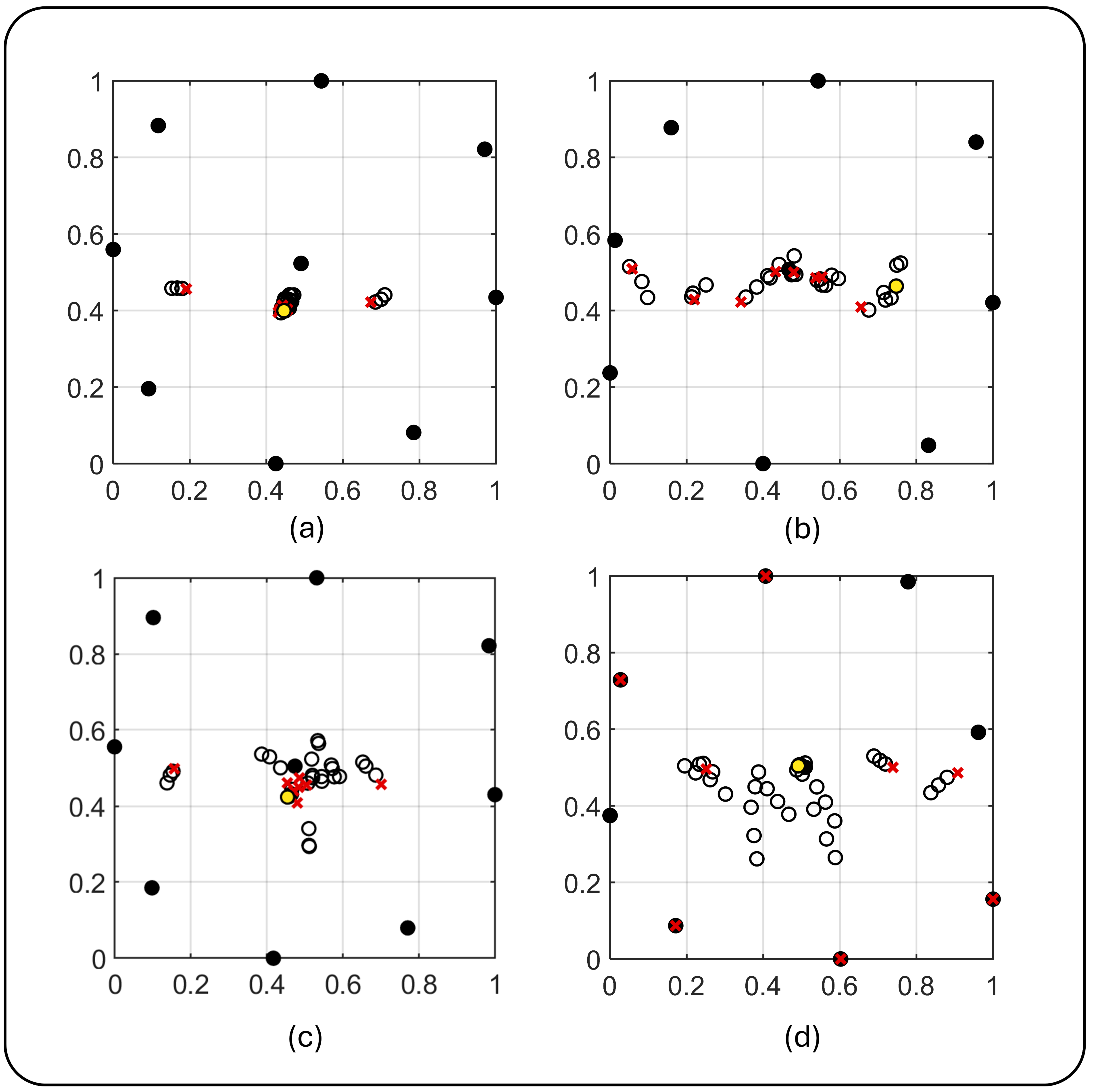}
\caption{Ablation study on gradient fixing strategies. Panels (a) and (c) represent the compliance minimization problem, while (b) and (d) represent the maximum stress minimization problem. The top row shows embedded samples without fixing the gradient, and the bottom row shows embedded samples with the gradient fixed. In all cases, $\eta_{\max}$ is reduced to 0.01 after the 5th iteration.
It can be observed that without fixing the gradient, many sampled points tend to cluster or overlap, leading to a reduction in the number of effective samples. In contrast, fixing the gradient yields better-separated samples, ensuring a more diverse and effective sampling distribution. The corresponding numerical example is presented in the first and second column of Fig.~\ref{compare1}.}
\label{fix}
\end{figure}

\paragraph{Ablation on Gradient Fixing Strategy} To validate the effectiveness of our gradient fixing strategy, we conducted an ablation study comparing the embedding spaces with and without this mechanism. As shown in Fig.~\ref{fix}, Panels (a) and (c) correspond to the compliance minimization problem, while (b) and (d) represent maximum stress minimization. The results in the top row indicate that without fixing the gradient, the sampled points tend to cluster or overlap significantly, leading to a loss of effective sample diversity. In contrast, the bottom row demonstrates that fixing the gradient facilitates better separation among samples. This strategy ensures a more diverse and effective sampling distribution, thereby preventing the algorithm from getting trapped in local clusters.

\section{Conclusion}
\label{sec:restrictions_outlook}

This work presents a framework for visualizing and quantifying the complexity of TO problems. We map high-dimensional design samples to a two dimensional embedding space and combine the embedding coordinates with objective values to form a visualization surface. Based on this surface, we quantify non-linearity using the index $\mathcal{I}_{\mathrm{NL}}$, where the convex envelope is approximated in by the lower convex hull.

Across structural, thermal, and flow benchmarks, we report extensive qualitative and quantitative results that compare visualization surfaces and analyze how parameter configurations and stabilization strategies influence non-linearity.

Several limitations remain. Projecting a high-dimensional space data to two dimensions inevitably discards information. The visualization surface depends on the sampling configuration, including the choice of initial points and the  $\eta_{\max}$. Moreover, the current study relies mainly on a single scalar index $\mathcal{I}_{\mathrm{NL}}$; while it captures overall non-linearity, it does not fully characterize the distribution of local minima or detailed curvature variations. Finally, our tests focus on steady problems at moderate scales; strongly non-linear material models and very large 3D systems remain to be examined.

\appendix
\section{Numerical Implementation Details}
\label{app1}

The topology optimization framework developed in this work is informed by well-established open-source codes for solid mechanics, fluid flow, and heat conduction problems \cite{code3DTO,heatobjectivefunctions,codeFlow}. The optimization algorithm employs the Method of Moving Asymptotes (MMA), utilizing the original code provided by Svanberg \cite{MMA1987}, with further implementation details available in \cite{svanberg2007mma}.








% ...
\bibliographystyle{elsarticle-num-names}     
\bibliography{ltexpprt_references}           

@article{eariest_TO,
  title={Generating optimal topologies in structural design using a homogenization method},
  author={Bends{\o}e, Martin Philip and Kikuchi, Noboru},
  journal={Computer methods in applied mechanics and engineering},
  volume={71},
  number={2},
  pages={197--224},
  year={1988},
  publisher={Elsevier}
}

@book{redbook,
  author    = {Bends{\o}e, Martin Philip and Sigmund, Ole},
  title     = {Topology Optimization: Theory, Methods, and Applications},
  edition   = {2nd},
  year      = {2004},
  publisher = {Springer},
  address   = {Berlin, Heidelberg},
  doi       = {10.1007/978-3-662-05086-6},
  url       = {https://doi.org/10.1007/978-3-662-05086-6},
  isbn      = {978-3-662-05086-6} 
}

@article{ole_review,
  title={Topology optimization approaches: A comparative review},
  author={Sigmund, Ole and Maute, Kurt},
  journal={Structural and multidisciplinary optimization},
  volume={48},
  number={6},
  pages={1031--1055},
  year={2013},
  publisher={Springer}
}

@article{ole_simp,
  title={Optimal shape design as a material distribution problem},
  author={Bends{\o}e, Martin P},
  journal={Structural optimization},
  volume={1},
  number={4},
  pages={193--202},
  year={1989},
  publisher={Springer}
}

@article{namesimp,
  title={The COC algorithm, Part II: Topological, geometrical and generalized shape optimization},
  author={Zhou, Ming and Rozvany, George IN},
  journal={Computer methods in applied mechanics and engineering},
  volume={89},
  number={1-3},
  pages={309--336},
  year={1991},
  publisher={Elsevier}
}

@article{stressconvergeslow,
  title={Stress-constrained topology optimization based on maximum stress measures},
  author={Yang, Dixiong and Liu, Hongliang and Zhang, Weisheng and Li, Shi},
  journal={Computers \& Structures},
  volume={198},
  pages={23--39},
  year={2018},
  publisher={Elsevier}
}

@article{3field2022,
  title={Three-field floating projection topology optimization of continuum structures},
  author={Huang, Xiaodong and Li, Weibai},
  journal={Computer Methods in Applied Mechanics and Engineering},
  volume={399},
  pages={115444},
  year={2022},
  publisher={Elsevier}
}

@article{3field2024,
  title={Reformulation for stress topology optimization of continuum structures by floating projection},
  author={Huang, Xiaodong and Li, Weibai and Nabaki, Khodamorad and Yan, Xiaolei},
  journal={Computer Methods in Applied Mechanics and Engineering},
  volume={423},
  pages={116870},
  year={2024},
  publisher={Elsevier}
}

@article{3field2025,
  title={Automatic projection parameter increase for three-field density-based topology optimization},
  author={Dunning, Peter and Wein, Fabian},
  journal={Structural and multidisciplinary optimization},
  volume={68},
  number={2},
  pages={1--19},
  year={2025},
  publisher={Springer}
}

@article{qualitymetricsvisual,
  title={Quality metrics in high-dimensional data visualization: An overview and systematization},
  author={Bertini, Enrico and Tatu, Andrada and Keim, Daniel},
  journal={IEEE Transactions on Visualization and Computer Graphics},
  volume={17},
  number={12},
  pages={2203--2212},
  year={2011},
  publisher={IEEE}
}

@article{goodfellowlandscape,
  title={Qualitatively characterizing neural network optimization problems},
  author={Goodfellow, Ian J and Vinyals, Oriol and Saxe, Andrew M},
  journal={arXiv preprint arXiv:1412.6544},
  year={2014}
}

@article{landscapevisualizing,
  title={Visualizing the loss landscape of neural nets},
  author={Li, Hao and Xu, Zheng and Taylor, Gavin and Studer, Christoph and Goldstein, Tom},
  journal={Advances in neural information processing systems},
  volume={31},
  year={2018}
}

@article{olecostwithmesh,
  title={On the usefulness of non-gradient approaches in topology optimization},
  author={Sigmund, Ole},
  journal={Structural and Multidisciplinary Optimization},
  volume={43},
  number={5},
  pages={589--596},
  year={2011},
  publisher={Springer}
}

@article{olegooddesignscale,
  title={Numerical instabilities in topology optimization: a survey on procedures dealing with checkerboards, mesh-dependencies and local minima},
  author={Sigmund, Ole and Petersson, Joakim},
  journal={Structural optimization},
  volume={16},
  number={1},
  pages={68--75},
  year={1998},
  publisher={Springer}
}

@article{LFT1,
  title={Legendre-Fenchel transforms in a nutshell},
  author={Touchette, Hugo},
  journal={URL http://www. maths. qmul. ac. uk/ht/archive/lfth2. pdf},
  pages={25},
  year={2005}
}

@article{fieldTO,
  title={A detailed introduction to density-based topology optimisation of fluid flow problems with implementation in MATLAB},
  author={Alexandersen, Joe},
  journal={Structural and Multidisciplinary Optimization},
  volume={66},
  number={1},
  pages={12},
  year={2023},
  publisher={Springer}
}

@article{ole_filter,
  title={Filters in topology optimization based on Helmholtz-type differential equations},
  author={Lazarov, Boyan Stefanov and Sigmund, Ole},
  journal={International journal for numerical methods in engineering},
  volume={86},
  number={6},
  pages={765--781},
  year={2011},
  publisher={Wiley Online Library}
}

@article{heatobjectivefunctions,
  title={A study on practical objectives and constraints for heat conduction topology optimization},
  author={Lohan, Danny J and Dede, Ercan M and Allison, James T},
  journal={Structural and Multidisciplinary Optimization},
  volume={61},
  number={2},
  pages={475--489},
  year={2020},
  publisher={Springer}
}

@article{PCAacceleratingphase-field,
  title={Accelerating phase-field predictions via recurrent neural networks learning the microstructure evolution in latent space},
  author={Hu, C and Martin, S and Dingreville, R},
  journal={Computer Methods in Applied Mechanics and Engineering},
  volume={397},
  pages={115128},
  year={2022},
  publisher={Elsevier}
}

@article{reduceDiminverse,
  title={Stabilizing and solving unique continuation problems by parameterizing data and learning finite element solution operators},
  author={Burman, Erik and Larson, Mats G and Larsson, Karl and Lundholm, Carl},
  journal={Computer Methods in Applied Mechanics and Engineering},
  volume={444},
  pages={118111},
  year={2025},
  publisher={Elsevier}
}

@article{PCADDTD,
  title={Data-driven topology design based on principal component analysis for 3D structural design problems: J. Yang et al.},
  author={Yang, Jun and Yaji, Kentaro and Yamasaki, Shintaro},
  journal={Structural and Multidisciplinary Optimization},
  volume={68},
  number={5},
  pages={103},
  year={2025},
  publisher={Springer}
}

@article{LFT2,
  title={A linear-time approximate convex envelope algorithm using the double Legendre--Fenchel transform with application to phase separation},
  author={Contento, Lorenzo and Ern, Alexandre and Vermiglio, Rossana},
  journal={Computational Optimization and Applications},
  volume={60},
  number={1},
  pages={231--261},
  year={2015},
  publisher={Springer}
}

@article{densitystart,
  title={Optimal shape design as a material distribution problem},
  author={Bends{\o}e, Martin P},
  journal={Structural optimization},
  volume={1},
  number={4},
  pages={193--202},
  year={1989},
  publisher={Springer}
}

@article{dasistart,
  title={Topology optimization of fluids in Stokes flow},
  author={Borrvall, Thomas and Petersson, Joakim},
  journal={International journal for numerical methods in fluids},
  volume={41},
  number={1},
  pages={77--107},
  year={2003},
  publisher={Wiley Online Library}
}

@article{design_space_dimensionality_reduction,
  title={A survey on design-space dimensionality reduction methods for shape optimization},
  author={Serani, Andrea and Diez, Matteo},
  journal={Archives of Computational Methods in Engineering},
  pages={1--28},
  year={2025},
  publisher={Springer}
}

@article{svanberg2007mma,
  title={MMA and GCMMA-two methods for nonlinear optimization},
  author={Svanberg, Krister},
  journal={vol},
  volume={1},
  pages={1--15},
  year={2007}
}

@article{MMA1987,
  title={The method of moving asymptotes—a new method for structural optimization},
  author={Svanberg, Krister},
  journal={International journal for numerical methods in engineering},
  volume={24},
  number={2},
  pages={359--373},
  year={1987},
  publisher={Wiley Online Library}
}

@article{scatteredInterpolant,
  title={Scattered data interpolation methods for electronic imaging systems: a survey},
  author={Amidror, Isaac},
  journal={Journal of electronic imaging},
  volume={11},
  number={2},
  pages={157--176},
  year={2002},
  publisher={Society of Photo-Optical Instrumentation Engineers}
}

@article{convhulln,
  title={The quickhull algorithm for convex hulls},
  author={Barber, C Bradford and Dobkin, David P and Huhdanpaa, Hannu},
  journal={ACM Transactions on Mathematical Software (TOMS)},
  volume={22},
  number={4},
  pages={469--483},
  year={1996},
  publisher={Acm New York, NY, USA}
}

@article{code3DTO,
  title={An efficient 3D topology optimization code written in Matlab},
  author={Liu, Kai and Tovar, Andr{\'e}s},
  journal={Structural and multidisciplinary optimization},
  volume={50},
  number={6},
  pages={1175--1196},
  year={2014},
  publisher={Springer}
}

@article{codeFlow,
  title={A detailed introduction to density-based topology optimisation of fluid flow problems with implementation in MATLAB},
  author={Alexandersen, Joe},
  journal={Structural and Multidisciplinary Optimization},
  volume={66},
  number={1},
  pages={12},
  year={2023},
  publisher={Springer}
}

@article{mds,
  title={Multidimensional scaling by optimizing goodness of fit to a nonmetric hypothesis},
  author={Kruskal, Joseph B},
  journal={Psychometrika},
  volume={29},
  number={1},
  pages={1--27},
  year={1964},
  publisher={Springer-Verlag}
}

@article{pca,
  title={Principal component analysis},
  author={Kim, KI and Jung, K and Kim, HJ},
  journal={Signal Processing},
  volume={9},
  number={2},
  pages={40--42},
  year={2002}
}

@article{tsne,
  title={Visualizing data using t-SNE},
  author={Maaten, Laurens van der and Hinton, Geoffrey},
  journal={Journal of machine learning research},
  volume={9},
  number={Nov},
  pages={2579--2605},
  year={2008}
}
\end{document}